\documentclass[12pt]{article}

\usepackage{makeidx}     
\usepackage{graphicx}    
\graphicspath{ {Figures/}, {Figures/Part1/}, {Figures/Part2/} }
\usepackage{xcolor}
\usepackage{wrapfig}
\usepackage{subcaption}
\usepackage{url}
\usepackage{multicol}    
\usepackage{amsmath}
\usepackage{amssymb}
\usepackage{alltt}
\usepackage{afterpage}
\usepackage{nicefrac}
\usepackage{url}
\usepackage[colorlinks, 
            pdfpagelabels, 
            pdfstartview = FitH, 
            bookmarksopen = true, 
            bookmarksnumbered = true,
            linkcolor = black,
            plainpages = false,
            hypertexnames = false,
            citecolor = black,
            urlcolor = black]{hyperref}
\usepackage{tikz}
\usepackage{pgfplots}
\usepackage{pgfplotstable}
\pgfplotsset{compat=1.10}
\usepgfplotslibrary{fillbetween}
\usetikzlibrary{calc,
                fadings, 
                intersections,
                circuits.ee.IEC}

\newcommand{\typenot}[1]{\typeout{was auf Seite \thepage nicht gedruckt}}

\newcommand{\D}{{\rm d}}
\newcommand{\ac}{\text{AC}\left( [0, T] \right)}
\newcommand{\bv}{\text{SBV}\left( [0, T] \right)}

\newtheorem{theorem}{Theorem}
\DeclareMathOperator{\im}{im}
\DeclareMathOperator{\spann}{span}
\DeclareMathOperator{\proj}{prox}
\DeclareMathOperator{\dist}{dist}
\newcommand{\oA}{{\cal A}}
\newcommand{\oB}{{\cal B}}
\newcommand{\zero}{\mbox{\boldmath $0$}}

\newcommand{\Id}{\mbox{\boldmath $I$}}

\newcommand{\dt}{\, \D t}

\newcommand{\wB}{\mbox{\boldmath $w$}}

\newcommand{\ibv}{\mbox{\boldmath $\eta$}}

\newcommand{\MM}{\mbox{\boldmath $M$}}
\newcommand{\ff}{\mbox{\boldmath $f$}}

\newcommand{\qe}{\mbox{\boldmath $q$}}

\newcommand{\Po}{\mbox{\boldmath $ q$}}
\newcommand{\Poc}{q}

\newcommand{\Vo}{\mbox{\boldmath $ v$}}

\newcommand{\Imp}{\mbox{\boldmath $ p$}}

\newcommand{\Mr}{\mbox{\boldmath $ M$}}
\newcommand{\fr}{\mbox{\boldmath $ f$}}

\newcommand{\GG}{\mbox{\boldmath $ G$}}
\newcommand{\gG}{\mbox{\boldmath $ g$}}

\newcommand{\Gr}{\mbox{\boldmath $ G$}}
\newcommand{\gr}{\mbox{\boldmath $ g$}}
\newcommand{\zr}{\mbox{\boldmath $ \kappa$}}

\newcommand{\La}{\mbox{\boldmath $ \lambda$}}
\newcommand{\Lam}{\mbox{\boldmath $ \mu$}}

\newcommand{\Zs}{\mbox{\boldmath $ q$}}

\newcommand{\Lac}{\lambda}

\newcommand{\pidy}{\mbox{\boldmath $ \delta$}}

\newcommand{\pinb}{\mbox{\boldmath $ \theta$}}

\newcommand{\svec}[1]{\mbox{\boldmath $#1$}}
\newcommand{\Vs}{\mbox{\boldmath $ v$}}

\newcommand{\hh}{\tau}

\newcommand{\xa}{\mbox{\boldmath $ x$}}
\newcommand{\xac}{x}
\newcommand{\Ea}{\mbox{\boldmath $ E$}}
\newcommand{\Ba}{\mbox{\boldmath $ H$}}
\newcommand{\ca}{\mbox{\boldmath $ c$}}
\newcommand{\Ua}{\mbox{\boldmath $ U$}}
\newcommand{\Va}{\mbox{\boldmath $ V$}}
\newcommand{\Ca}{\mbox{\boldmath $ C$}}

\newcommand{\Na}{\mbox{\boldmath $ N$}}
\newcommand{\Fa}{\mbox{\boldmath $ \phi$}}

\newcommand{\ya}{\mbox{\boldmath $ y$}}
\newcommand{\za}{\mbox{\boldmath $ z$}}
\newcommand{\yac}{y}
\newcommand{\zac}{z}
\newcommand{\yav}{\mbox{\boldmath $ v$}}
\newcommand{\fa}{\mbox{\boldmath $ a$}}
\newcommand{\ga}{\mbox{\boldmath $ b$}}
\newcommand{\Ga}{\mbox{\boldmath $ B$}}

\newcommand{\FFa}{\mbox{\boldmath $ F$}}
\newcommand{\Par}{\mbox{\boldmath $ \psi$}}
\newcommand{\DPar}{\mbox{\boldmath $ \Psi$}}
\newcommand{\xlo}{\mbox{\boldmath $ \xi$}}

\newcommand{\Xa}{\mbox{\boldmath $ X$}}

\newcommand{\Ya}{\mbox{\boldmath $ Y$}}
\newcommand{\Za}{\mbox{\boldmath $ Z$}}

\begin{document}

\begin{center}
{\Large \bf Differential-Algebraic
 Equations  and Beyond: \\[4mm] }
{ \bf From Smooth to Nonsmooth Constrained Dynamical Systems \\[4mm] }
Jan Kleinert \\
{\small Hochschule Bonn-Rhein-Sieg} \\
{\small Grantham-Allee 20, D-53757 Sankt Augustin} \\
{\small \tt jan.kleinert@h-brs.de}\\[2mm]
{\small and} \\[2mm]
{\small German Aerospace Center}\\
{\small Simulation and Software Technology} \\
{\small Linder H{\"o}he, D-51147 Cologne} \\
{\small \tt jan.kleinert@dlr.de}\\[2mm]
Bernd Simeon \\ 
{\small Felix-Klein-Zentrum, TU Kaiserslautern }\\
{\small D-67663 Kaiserslautern, Germany }\\
{\small \tt simeon@mathematik.uni-kl.de} \\[2mm]
\today \\
\mbox{\quad}
\end{center}
\begin{quote}
{\bf Abstract:\/} The present article 
presents a summarizing view at differential-algebraic equations (DAEs)
and analyzes how new application fields and corresponding mathematical models lead to innovations both in theory and in numerical analysis for
this problem class. Recent numerical methods for nonsmooth dynamical systems subject to unilateral contact and friction illustrate the topicality of this 
development. 
\end{quote}

\medskip


\section{Introduction}
Differential-Algebraic Equations (DAEs) \index{differential-algebraic equation}
are a prominent example for applica\-tion-driven research 
that leads to new concepts and methodology in mathematics.
Until the early 1980s of the last century, this topic was widely unknown but the introduction of 
powerful simulation software in  electrical and mechanical engineering created a strong demand for the analysis and numerical solution of dynamical systems with constraints. 
Mathematicians all over the world then started to work on DAEs, which resulted in an avalanche of research over the following decades. Meanwhile, most issues have been resolved and sophisticated numerical methods have been found, but despite this maturity, the field is constantly expanding due to the ongoing trend to model complex phenomena in science and engineering by means of
differential equations and additional constraints that stem from network structures, boundary and coupling conditions or physical conservation properties.

In  this paper, we strive for a survey of this development and even more also 
discuss recent extensions towards dynamical systems that are non\-smooth. The latter topic is quite timely and arises, e.g., in the modeling of granular material.
Clearly, our approach is exemplary in nature and does not aim at completeness. There is, moreover,
a strong bias on mechanical systems.
 Those readers who would like to know more about the topic of DAEs and the rich oeuvre that has accumulated over the years are referred to 
the monographs of Brenan, Campbell \& Petzold \cite{BrCP96}, 
Griepentrog \& M\"arz \cite{GrMa86},
Hairer \& Wanner \cite{HaWa96}, 
Kunkel \& Mehrmann \cite{KuMe06}, 
Lamour, M\"arz \& Tischendorf \cite{LaMaTi13},
and to the survey of Rabier \& Rheinboldt \cite{RaRhei02}. Refer to the textbooks by Acary and Brogliato \cite{acary2008,brogliato2016} and references therein for a more elaborate literature survey on nonsmooth dynamical systems.
 
The paper is organized as follows. Section~\ref{sect:dae} outlines the emergence of DAEs and shows how it was influenced by applications from electrical and mechanical engineering problems. Section~\ref{sect:dae-results} then highlights some major results and numerical methods that were developed in the DAE context. Some extensions to classical DAEs for partial and stochastic differential equations are discussed in Section~\ref{sect:beyond}.
Section~\ref{sect:nonsmoothdynsys} discusses nonsmooth dynamical systems. An overview of the most important theoretical concepts is given in Section~\ref{sect:dviHistory}, the equations of motion for a unilaterally constrained mechanical system are motivated in Section~\ref{sect:nonsmoothMechanics} and finally, Section~\ref{sect:nonsmoothNumerics} provides some insights into numerical methods and recent developments in the field of nonsmooth dynamical systems.
%
%
\section{Differential-algebraic equations} \label{sect:dae}
In this section, we describe how the DAEs became a hot topic in the 1980s
and 1990s and then go further back in time to the works of Kirchhoff \cite{Kirchhoff1847}
and Lagrange  \cite{Lag88} who introduced differential equations with constraints in order to model electric circuits and mechanical systems.
\subsection{How the topic of DAEs emerged}
In the beginning of the 1980s, the term `DAE' was widely unknown in mathematics. But this changed rapidly, due to an increasing demand in several engineering fields but also due to  the pioneering work of Bill Gear.
The first occurrence of the term \emph{Differential-Algebraic Equation} 
\index{differential-algebraic equation} 
can be found in the title of Gear's paper  \emph{Simultaneous numerical solution of differential-algebraic equations} \cite{Gear71} from 1971, and
in the same year his famous book 
\emph{Numerical Initial Value Problems in Ordinary Differential Equations}
\cite{Gear71a} appeared where he already considers examples from electric circuit analysis  in the form 
\begin{equation}\label{dae:li}
   \Ea \, \dot{\xa} = \Fa(\xa,t)
\end{equation}
with possibly singular capacitance matrix $\Ea \in \mathbb{R}^{n_\xac \times n_\xac}$
and right hand side function~$\Fa$. For a regular matrix $\Ea$, (\ref{dae:li}) can be 
easily converted into a system of ordinary differential equations (ODEs). Otherwise, however,
a DAE arises that calls for new approaches, both in theory and in numerical analysis.


\begin{figure}[t]
\center
\includegraphics[height=4.3cm]{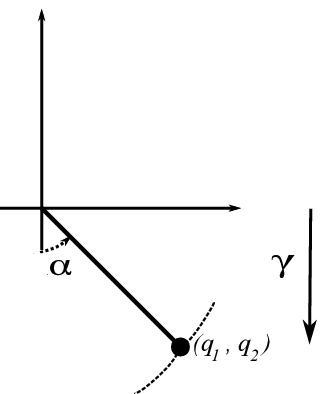}\quad
\includegraphics[height=5cm]{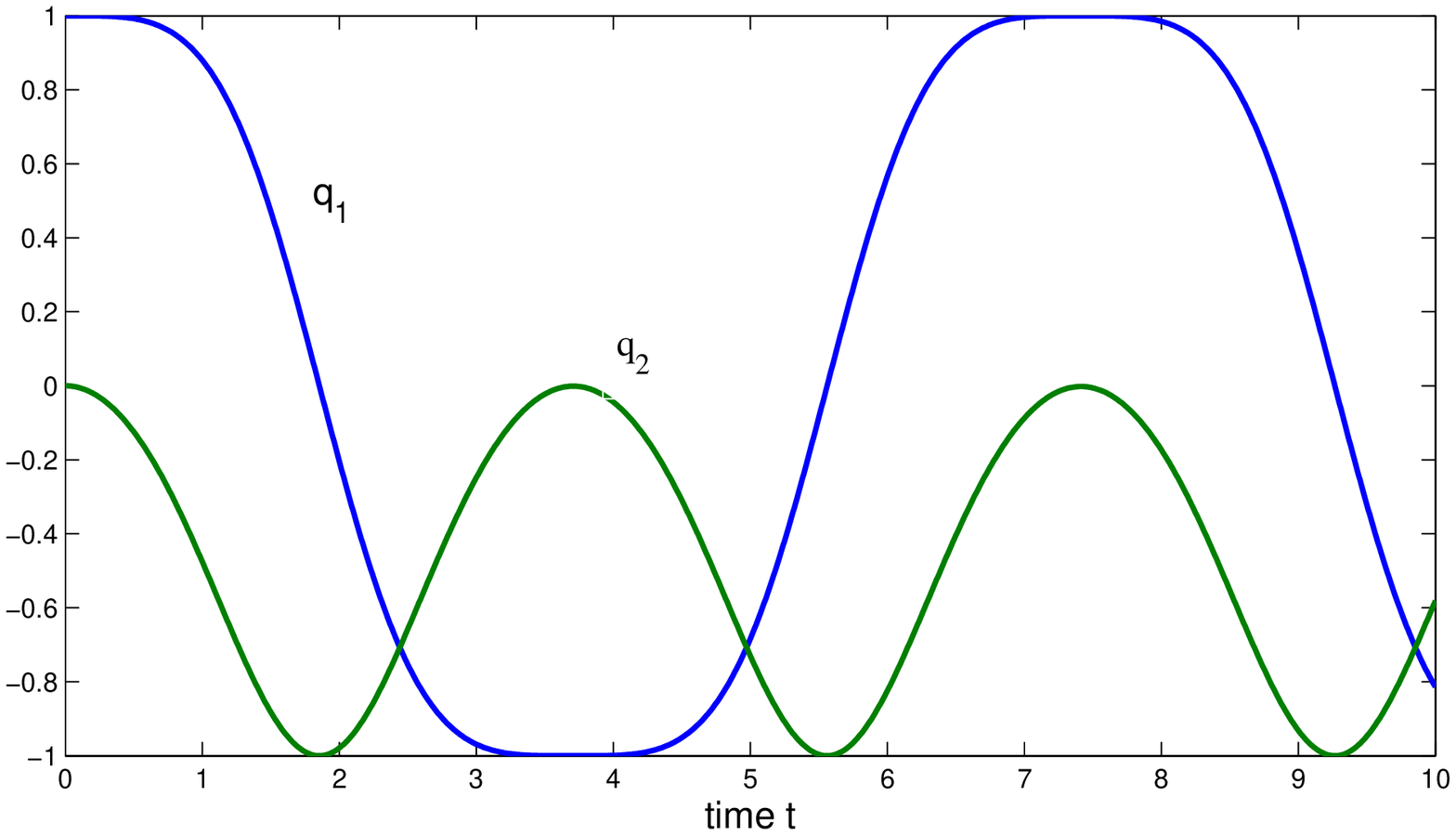}
   \caption{The mathematical pendulum. 
} \label{fig:pendulum}
\end{figure}

Two application fields, namely electric circuit analysis and constrained mechanical systems, 
have been the major driving forces for the development of DAEs. Below, this statement will be made more explicit by looking at the corresponding modeling concepts. 
Bill Gear had the farsightedness to perceive  very early the importance of these modeling approaches for today's simulation software. During an Oberwolfach workshop in 1981, he suggested to 
study the \emph{mathematical pendulum in Cartesian coordinates}
\begin{subequations}\label{mathpend}
\begin{eqnarray}
   \ddot{\Poc}_1 & = & - 2 \Poc_1 \Lac, \label{mathpenda} \\
   \ddot{\Poc}_2 & = & - \gamma - 2 \Poc_2 \Lac, \label{mathpendb} \\
            0           & = & \Poc_1^2 + \Poc_2^2 -1 \label{mathpendc} 
\end{eqnarray}
\end{subequations}
that describes the motion of a mass point with coordinates $(\Poc_1,\Poc_2)$ in the plane subject to a constraint. The constraint models the massless rod of length $1$ that connects the mass point to a pivot placed in the origin of the coordinate system, Fig.~\ref{fig:pendulum}. The motion of the mass point is then determined by the gravity (parameter $\gamma$) and by the constraint forces that are expressed in terms of the unknown Lagrange multiplier $\lambda$.

The DAE (\ref{mathpend}) is an example for the Lagrange equations of the first kind that we will discuss below. By introducing velocity variables, it can be easily converted to a system of first order that fits into the class of linear-implicit DAEs (\ref{dae:li}).
The most general form of a DAE is a
\emph{fully implicit system}
\begin{equation}\label{dae:fi}
         \FFa(\dot{\xa},\xa,t) = \zero
\end{equation}
with state variables $\xa(t) \in \mathbb{R}^{n_\xac}$ and 
a nonlinear, vector-valued function $\FFa$ of corresponding dimension.
Clearly, if the $n_\xac \times n_\xac$ Jacobian $\partial \FFa / \partial \dot{\xa}$
is invertible, then by the implicit function theorem, it is theoretically possible to transform (\ref{dae:fi}), at least locally, to an explicit system of ODEs.
If $\partial \FFa / \partial \dot{\xa}$ is singular, however, (\ref{dae:fi}) 
constitutes a DAE.

Linda Petzold, a student of Bill Gear, continued and extended his pioneering work in various
directions. In particular, the development of the DASSL code (the `Differential-Algebraic System SoLver')  that she had started in the early 1980s \cite{Petz82,BrCP96} set a corner stone that still persists today. DASSL is based on the
Backward Differentiation Formulas (BDF), which are also popular for solving systems of stiff ODEs.
The extension of the BDF methods to implicit systems (\ref{dae:fi})  is intriguingly simple. 
One replaces the differential operator~$\D/\D t$ in (\ref{dae:fi})  by the difference operator
\index{BDF methods}
\begin{equation}\label{dae:bdfk}
\varrho \xa_{n+k} :=
\sum_{i=0}^k \alpha_i \xa_{n+i} = \hh \dot{\xa}(t_{n+k}) 
+ {\cal O}(\hh^{k+1}) 
\end{equation} 
where $\xa_{n+i}$ stands for the discrete approximation of $\xa(t_{n+i})$ with stepsize  $\hh$ and where the $\alpha_{i}, \, i=0,\ldots,k$, denote the method coefficients. 
Using the finite difference approximation $\varrho \xa_{n+k}/\hh$ of the time derivative,
the numerical solution of the DAE (\ref{dae:fi}) 
\index{differential-algebraic equation!implicit} 
then boils down to solving the nonlinear system
\begin{equation}\label{dae:finum}
         \FFa\left(\frac{\varrho\xa_{n+k}}{\hh},\xa_{n+k},t_{n+k}\right) = \zero
\end{equation}
for $\xa_{n+k}$ in each time step.
This is exactly the underlying idea of the DASSL code, see  Fig.~\ref{fig:dassl}.

Soon after the first release of the DASSL code, it became very popular among 
engineers and mathematicians. For quite some problems, however, the code 
would fail, which in turn triggered new research in numerical analysis in order to
understand such phenomena. As it turned out, the notion of an {\em index} 
of the DAE (\ref{dae:fi}) was the key to obtain further insight.

\begin{figure}[t]
{\small \begin{verbatim}
      SUBROUTINE DDASSL (RES,NEQ,T,Y,YPRIME,TOUT,INFO,RTOL,ATOL,
     +     IDID,RWORK,LRW,IWORK,LIW,RPAR,IPAR,JAC)
C***BEGIN PROLOGUE  DDASSL
C***PURPOSE  This code solves a system of differential/algebraic
C            equations of the form G(T,Y,YPRIME) = 0.
\end{verbatim}
}
   \caption{Calling sequence of the DASSL code \cite{Petz82,BrCP96} that has
had an enormous impact on the subject of DAEs
\index{differential-algebraic equation} 
 and that is still in wide use today. 
The original code is written in FORTRAN77 in double precision. A recent implementation in C is part of the SUNDIALS suite of codes \cite{HBetal05}.
} \label{fig:dassl}
\end{figure}

Gear \cite{Gear88,Gear90} introduced what we call today the 
{\em differentiation index}%
\index{index!differential}. 
This non-negative integer $k$ is defined by 
\begin{description}
\item[$k = 0$:] If $\partial \FFa /\partial \dot{\xa}$ is non-singular, the index is $0$.
\item[$k > 0$:] Otherwise, consider the system of equations
\index{differential-algebraic equation!implicit}
\begin{eqnarray}
		   \FFa(\dot{\xa},\xa,t) & = & \zero, \nonumber \\
		     \frac{ \D}{ \D t} \FFa(\dot{\xa},\xa,t)
		      = \frac{\partial}{\partial \dot{\xa}} \FFa(\dot{\xa},\xa,t)\, \xa^{(2)} + \ldots
		      & =& \zero, \label{dae:di} \\
				 & \vdots   &                 \nonumber      \\
		     \frac{ \D^s}{ \D t^s} \FFa(\dot{\xa},\xa,t)
		      = \frac{\partial}{\partial \dot{\xa}} \FFa(\dot{\xa},\xa,t)\, \xa^{(s+1)} + \ldots
		      & =& \zero \nonumber
\end{eqnarray}
as a system in the separate dependent
variables $\dot{\xa}, \xa^{(2)}, \ldots, \xa^{(s+1)}$,
with $\xa$ and $t$  as independent variables. Then the index $k$ is the smallest
$s$ for which it is possible, using algebraic manipulations only, to extract 
an ordinary differential equation $\dot{\xa} = \svec{\psi}(\xa,t)$ 
(the underlying ODE)
from~(\ref{dae:di}).
\end{description}

Meanwhile other notions of an index have emerged, but despite its ambiguity 
with respect to the algebraic manipulations, the differentiation 
index is still the most popular and widespread tool to classify DAEs.

In the next chapter, other index concepts and their relation to the differential index will be addressed, and also more protagonists will enter the stage. 
This first section on the early days of DAEs closes now with a look at the application fields that set the ball rolling.

\subsection{Electric circuits}

In 1847, Kirchhoff first published his \emph{circuit laws} that describe the 
conservation properties of electric circuits \cite{Kirchhoff1847}. 
These laws consist of the current law and the voltage law, which both follow from Maxwell's equations of electro-dynamics. When these laws are applied to circuits with time-dependent behavior, the corresponding equations are typically given as a linear-implicit system (\ref{dae:li}).
Often, the structure even turns out to be a linear constant coefficient DAE 
\index{differential-algebraic equation!linear constant coefficient} 
\begin{equation}\label{dae:constc}
 \Ea \dot{\xa} + \Ba \xa = \ca(t) 
\end{equation}
with matrices $\Ea,\Ba \in \mathbb{R}^{n_\xac \times n_\xac}$
and a time-dependent source term~$\ca(t) \in \mathbb{R}^{n_\xac}$.

\begin{figure}[t]
\center
\includegraphics[height=4cm]{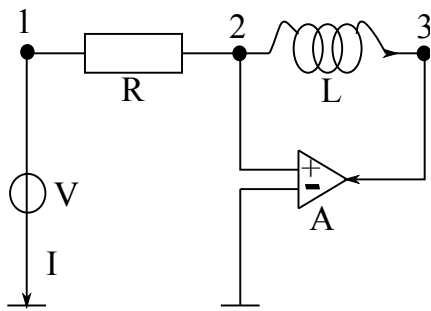}
   \caption{Differentiator circuit. 
} \label{fig:different}
\end{figure}

An example for such an electric circuit is the {\em differentiator\/}  \cite{GHR2000} shown in
Fig.~\ref{fig:different}. It consists of a resistance $R$, an inductance $L$, an ideal operational amplifier $A=\infty$, and a given voltage source $V(t)$. The $n_\xac = 6$ unknowns read here $\xa = ( V_1, V_2, V_3, I, I_L, I_V)$ with voltages $V_i$ and currents $I, I_L, I_V$. From Kirchhoff's laws and the properties of the amplifier and the inductance one obtains the relations
\begin{eqnarray*}
    I + (V_1-V_2)/R & = & 0, \\
    -(V_1-V_2)/R + I_L & = & 0, \\
  -I_L + I_V & = & 0, \\
            V_1 & = & V(t), \\
            V_2 & = & 0, \\
  V_2 - V_3 & = & L \cdot \dot{I}_L .
\end{eqnarray*}
This linear system has the form 
(\ref{dae:constc}) with singular inductance matrix
\begin{equation}\label{opamp}
\Ea = 
\left( \begin{array}{ccccc}
 0 & 0 & \ldots & 0 & 0 \\
\vdots & \vdots & \ddots & \vdots & \vdots \\
0 & 0 & \ldots & 0 &  0 \\
0 & 0&0 & L & 0 
\end{array}\right) .
\end{equation}

If the matrix $\Ea$ was regular, it could be 
brought to the right hand side by formal inversion, ending up in a system of ODEs. 
Here however, it is singular and thus we
face a DAE problem.
 
Weierstrass and Kronecker were in Berlin at the same time as Kirchhoff, and
it is quite obvious to suppose that they knew his work \footnote{The relation of the work of Weierstrass and Kronecker to Kirchhoff's circuit laws was pointed out to me by Volker Mehrmann when we met in September 2014 during a Summer School on DAEs in Elgersburg, Germany.}.
 Weierstrass and later Kronecker were thus 
inspired to study such singular systems and provided
an elegant theory that is still fundamental today in order to
understand the specific properties of DAEs.

We assume that the {\em matrix pencil\/} $\mu \Ea + \Ba \in 
\mathbb{R}^{{n_x}\times{n_x}}[\mu]$ is regular. I.e., 
there exists $\mu \in \mathbb{C}$ such that the matrix
$\mu \Ea + \Ba$ is regular. Otherwise, the pencil is singular, and 
(\ref{dae:constc}) has either no or infinitely many solutions. This latter case 
has been first studied by Kronecker \cite{Kron90}, see also
\cite{Gant59,Camp82}.

If  $\mu \Ea + \Ba$ is regular, there exist nonsingular
matrices $\Ua$ and $\Va$ such that 
\begin{equation}\label{dae:weierstrass}
         \Ua\Ea\Va = \left( \begin{array}{cc} \Id & \zero \\
                                                  \zero & \Na
                                \end{array} \right), \quad
         \Ua\Ba\Va = \left( \begin{array}{cc} \Ca & \zero \\
                                                  \zero & \Id
                                \end{array} \right)                    
\end{equation}
where $\Na$ is a nilpotent matrix, $\Id$ the identity matrix, and
$\Ca$ a matrix that can be assumed to be in Jordan canonical form.
Note that the dimensions of these square blocks in (\ref{dae:weierstrass}) are uniquely determined.
The transformation 
(\ref{dae:weierstrass}) is called
the {\em Weierstrass canonical form\/} \cite{Weierstrass1868}.
It is a generalization of the 
Jordan canonical form and contains the essential structure of the
linear system (\ref{dae:constc}). 

In the Weierstrass canonical form (\ref{dae:weierstrass}), the 
\index{Weierstrass canonical form}
singularity of the DAE is represented by the nilpotent matrix
$\Na$. Its degree of nilpotency, i.e., the smallest positive integer
$k$ such that $\Na^k = \zero$,
plays a key role when studying closed-form solutions of the
linear system (\ref{dae:constc}) and is identical to 
the differentiation index of (\ref{dae:constc}).
\index{index!nilpotency} 

To construct a solution of (\ref{dae:constc}), we introduce new variables and right hand side vectors
\begin{equation}\label{dae:trafod}
       \Va^{-1} \xa =: \left( \begin{array}{c} \ya \\ \za 
                              \end{array} \right), \quad
        \Ua \ca =: \left( \begin{array}{c} \pidy \\ \pinb 
                              \end{array} \right).
\end{equation}
Premultiplying (\ref{dae:constc}) by $\Ua$ then leads to
the {\em decoupled system}
\begin{subequations}
\begin{eqnarray}
             \dot{\ya} + \Ca \ya & = & \pidy\,, \label{dae:deca} \\
             \Na \dot{\za} + \za  & = & \pinb \,. \label{dae:decb}
\end{eqnarray}
\end{subequations}
While the solution of the ODE (\ref{dae:deca}) follows by integrating
and results in an expression based on the matrix exponential
$\exp (-\Ca(t-t_0))$, the equation (\ref{dae:decb}) for $\za$
can be solved recursively by differentiating. More precisely,
it holds
\[   \Na \ddot{\za} + \dot{\za} = \dot{\pinb}
      \quad \Rightarrow \quad \Na^2  \ddot{\za} =  
            - \Na \dot{\za} + \Na \dot{\pinb} = \za -  \pinb + \Na
            \dot{\pinb} \,.
\]
Repeating the differentiation and multiplication by $\Na$,
we can eventually exploit the nilpotency and get 
\[   \zero = \Na^k \za^{(k)} = (-1)^k \za + \sum_{\ell = 0}^{k-1}
                               (-1)^{k-1-\ell} \Na^{\ell} \pinb^{(\ell)}.
\]
This implies the explicit representation
\begin{equation}\label{dae:solz}
           \za = \sum_{\ell = 0}^{k-1} (-1)^{\ell}
                               \Na^{\ell} \pinb^{(\ell)}.
\end{equation}

The above solution procedure illustrates several crucial points
about DAEs and how they differ from  ODEs. Remarkably, the linear constant coefficient case
displays already these points, and thus the work of Weierstrass and Kronecker represents still the fundament of DAE theory today.

We highlight two crucial points:

\begin{description}
\item[(i)] The solution of (\ref{dae:constc}) rests on $k-1$ 
  differentiation steps. This requires that the derivatives
  of certain components of $\pinb$  
  exist up to $\ell = k-1$. Furthermore, some components of $\za$ 
  may only be continuous but not differentiable depending
  on the smoothness of $\pinb$.
\item[(ii)] The components of $\za$ are directly given in terms
  of the right hand side data $\pinb$ and its derivatives. 
  Accordingly, the initial value $\za(t_0)= \za_0$ is fully determined
  by (\ref{dae:solz}) and, in contrast to $\ya_0$, cannot be chosen arbitrarily. 
  Initial values $(\ya_0,\za_0)$ where $\za_0$ satisfies
  (\ref{dae:solz}) are called {\em consistent.} The same terminology
  \index{consistent initial value}
  applies to the initial value $\xa_0$, which is consistent if, 
  after the transformation (\ref{dae:trafod}), $\za_0$ satisfies
  (\ref{dae:solz}).
\end{description}
Today, more than 150 years after the discoveries of Kirchhoff, electric circuit analysis remains one of the driving forces in the development of DAEs. 
The interplay of modeling and mathematical analysis is particularly important in this field, and the interested reader is referred to 
 G\"unther \& Feldmann \cite{GF99} 
and M\"arz \& Tischendorf \cite{MT97}
as basic works. The first simulation code that 
generated a model in differential-algebraic form was the SPICE package 
\cite{Nagel:M382}.

\begin{figure}[t]
\begin{center}
\includegraphics[width=10.5cm]{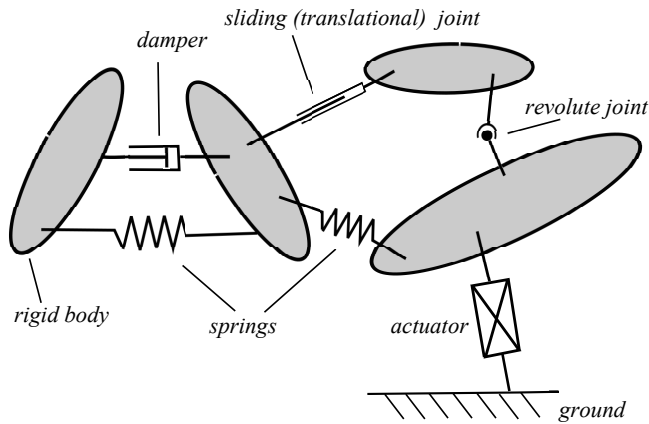}\mbox{\qquad}
\end{center}
   \caption{Sketch of a multibody system with rigid bodies and typical
   interconnections.
} \label{fig:mbsgeneral}
\end{figure}

\subsection{Constrained mechanical systems}

Even older than the DAEs arising from Kirchhoff's laws are the Euler-Lagrange equations. They were first published in Lagrange's famous work {\em M{\'e}canique analytique} \cite{Lag88} from 1788. 

Consider a mechanical system that consists of rigid bodies interacting
via springs, dampers, joints, and actuators, Fig.~\ref{fig:mbsgeneral}.
The bodies possess a certain geometry and mass while the interconnection elements are massless.
Let 
$\Po(t)\in \mathbb{R}^{n_\Poc}$ denote a vector that comprises 
the coordinates for position and orientation of all bodies in the system. 
Revolute, translational, universal, and spherical joints are examples for bondings in such a multibody system.
They may constrain the motion $\Po$ and hence determine its kinematics.
 
If constraints are present, we express the resulting conditions on 
$\Po$
in terms of $n_\Lac$ constraint
equations \index{constraint!at position level}
\begin{equation}\label{constrigid}
                \zero = \gr(\Po) \,.
\end{equation}
Obviously, a meaningful model requires $n_\Lac < n_\Poc$.
The equations (\ref{constrigid}) that restrict the motion 
$\Po$ are called  {\em holonomic constraints}, and the rectangular matrix
\[
    \Gr(\Po) := \frac{\partial \gr(\Po)}{\partial \Po} \,\, \in \mathbb{R}^{n_\Lac \times n_\Poc}
\]
is called the {\em constraint Jacobian}. 
\index{constraint!Jacobian}

Using both the redundant position variables $\Po$ and additional Lagrange multipliers 
$\La$ to describe the dynamics leads to the 
{\em equations of constrained mechanical motion}, also called the {\em Lagrange equations of the first kind\/} or the \index{Euler-Lagrange equations}
\emph{Euler-Lagrange equations}
\begin{subequations}\label{f1}
\begin{eqnarray}
            \Mr(\Po) \, \ddot{\Po}& =& \fr(\Po,\dot{\Po},t) - 
		      \Gr(\Po)^T \La  \,, \label{f1a} \\
             \zero & = & \gr(\Po)\,,   \label{f1b}
\end{eqnarray} 
\end{subequations}
where $\Mr(\Po)\in \mathbb{R}^{n_\Poc \times n_\Poc}$
stands for the {\em mass matrix} and \index{mass matrix} 
$\fr(\Po,\dot{\Po},t)\in \mathbb{R}^{n_\Poc}$ for the vector of {\em applied and internal forces}. 

The standard example for such a constrained mechanical system are the equations (\ref{mathpend}) of the mathematical pendulum. 
For a long time, it was common sense that the Euler-Lagrange equations should be transformed to the  {\em state space form},
also called  the {\em Lagrange equations of the second kind}. In case 
of the pendulum, this means that the Cartesian coordinates can be expressed as $\Poc_1 = \sin \alpha, \, \Poc_2 = - \cos \alpha$ with the angle $\alpha$ as
minimal coordinate, Fig.~\ref{fig:pendulum}. By inserting these relations into (\ref{mathpend}), the 
constraints and the Lagrange multiplier cancel, and one arrives at the second
order ODE
\begin{equation}\label{statespace}
   \ddot{\alpha} = - \gamma \sin \alpha 
\end{equation}
as state space form.

It seems obvious that
a state space form such as (\ref{statespace})  
constitutes a more appropriate and easier model than the differential-algebraic system (\ref{f1}),  or (\ref{mathpend}), respectively, in redundant coordinates. In practice, however, the state space form suffers from serious drawbacks:

The analytical complexity
of the constraint equations (\ref{constrigid}) makes it 
in various applications impossible to obtain a set of minimal coordinates that is
\index{coordinates!minimal} valid 
for all configurations of the multibody system.  Moreover, although we know from the theorem on implicit functions 
that such a set exists in a neighborhood of the current configuration, 
it might loose its validity when the configuration changes.
This holds in particular for 
multibody systems with so-called {\em closed kinematic loops}. 
\index{closed kinematic loop}

Even more, the modeling 
of subsystems
like electrical and hydraulic feedback controls, which are essential
for the performance of modern mechanical systems, is limited.
The differential-algebraic model, on the other hand, bypasses topological
analysis and
offers the choice of using a set of
coordinates $\Po$ that possess physical significance. 

This reasoning in favor of the differential-algebraic model (\ref{f1}) became more and more widespread in the 1980s, driven by the development of sophisticated software packages, so-called {\em multibody formalisms}. 
One of the first packages that fully exploited this new way of modeling is due to Haug \cite{Haug89}.
 
A look at the leading software tools in the field today shows a clear picture. Some of the codes
generate a differential-algebraic model whenever a constraint is present, while others try to generate 
a state space form as long as it is convenient. But 
the majority of the commercial products relies on the differential-algebraic 
approach as the most general way to handle complex technical applications \cite{GB94,Schi90}.
\index{multibody system}

The main difference between the DAEs arising from electric circuit analysis 
and the DAEs that model constrained mechanical systems is the richer structure of the latter.  
E.g., for {\em conservative multibody systems}, i.e., 
systems where the applied forces can be written as the gradient of a potential $U$, the Euler-Lagrange equations  (\ref{f1}) result from 
Hamilton's principle of least action
\begin{equation}\label{hamstarr}
  \int_{t_0}^{t_1} \Big( T - U - \gr(\Po)^T \La \Big) \dt 
  \rightarrow \mbox{ stationary !}
\end{equation}
where the kinetic energy possesses a representation as quadratic form
\index{kinetic energy}
\[  T(\Po,\dot{\Po}) = \frac{1}{2} \dot{\Po}^T \Mr(\Po) \dot{\Po} \,.
\]
In the least action principle (\ref{hamstarr}),
we observe the fundamental Lagrange multiplier technique \index{Lagrange!multiplier} for coupling constraints and dynamics \cite{Bri08}.
Extensions of the multiplier technique exist in various more general 
settings such as dissipative systems or even inequality constraints, see Section~\ref{sect:nonsmoothdynsys}. 

The pendulum equations (\ref{mathpend}) are the example for a constrained mechanical system. Though they simply describe the motion of a single mass point, several key properties of the Euler-Lagrange equations can already be  studied: the differential equations are of second order, the constraint equations are mostly nonlinear, and one observes a clear semi-explicit structure with differential variables $\Po$ and algebraic variables $\La$.

The Euler-Lagrange equations are of index $3$ and form the prototype for a
system of higher index. Index reduction techniques are thus required and in fact, already in 1972 this issue was addressed by Baumgarte 
\cite{Baum72}.  He observed that in (\ref{f1}), the Lagrange multipliers can be eliminated by differentiating the constraints twice. The first
differentiation leads to the constraints at velocity level
\begin{equation}\label{develo}
 \zero = \frac{\D}{\D t} \, \gr(\Po) = \Gr(\Po)\, \dot{\Po}. 
\end{equation}
A second differentiation step yields the {\em constraints at acceleration level}
\index{constraint!at acceleration level}
\begin{equation}\label{deaccel}
 \zero = \frac{\D^2}{\D t^2} \, \gr(\Po) = \Gr(\Po)\, \ddot{\Po}
          + \zr(\Po,\dot{\Po})\,, \quad \zr(\Po,\dot{\Po}) :=
          \frac{\partial \Gr(\Po)}{\partial \Po} (\dot{\Po},\dot{\Po})\,,
\end{equation}
where the two-form $\zr$ comprises additional derivative terms.
The combination of the dynamic equation 
\[
 \Mr(\Po) \, \ddot{\Po} = \fr(\Po,\dot{\Po},t) - 
		      \Gr(\Po)^T \La  
\]
with (\ref{deaccel}) results in a linear system for $\ddot{\Po}$ and
$\La$ with the saddle point matrix
\begin{equation}\label{vorausstarr}
\left( \begin{array}{cc} \Mr(\Po) & \Gr(\Po)^T \\ \Gr(\Po) & \zero
                        \end{array} \right) \in \mathbb{R}^{(n_\Poc+n_\Lac) \times
                              (n_\Poc+n_\Lac)}.
\end{equation}
For a well-defined multibody system, this matrix is invertible 
in a neighborhood of the solution, and in this way, the Lagrange multiplier 
can be computed as a function of $\Po$ and $\dot{\Po}$.

However, the well-known drift-off phenomenon requires additional stabilization measures, and Baumgarte came up with 
the idea  to combine 
original and differentiated constraints as 
\begin{equation}\label{baumgarte}
 \zero = \Gr(\Po)\, \ddot{\Po}
          + \zr(\Po,\dot{\Po}) + 2 \alpha \Gr(\Po) \dot{\Po} + \beta^2
\gr(\Po)
\end{equation}
with scalar parameters $\alpha$ and $\beta$.
The free parameters $\alpha$ and $\beta$ should be chosen in such
a way that
\begin{equation}\label{driftoffb}
   \zero = \ddot{\wB} + 2\alpha \dot{\wB} + \beta^2 \wB
\end{equation}
becomes an asymptotically stable equation, with 
$\wB(t) := \gr(\Po(t))$.

From today's perspective, the crucial point in Baumgarte's approach is the choice of 
the parameters. Nevertheless, it was the very beginning of a long series of works that tried to reformulate the Euler-Lagrange equations in such a way that the index is lowered while still maintaining the information of all
constraint equations.
For a detailed analysis of this stabilization 
and related techniques we refer to Ascher et al.~\cite{Asch95,Asch97}.

Another -- very early -- stabilization of the Euler-Lagrange equations is due to 
Gear, Gupta \& Leimkuhler \cite{GeGL85}. This formulation represents still the state-of-the-art in multibody dynamics.
It uses a formulation of the equations of motion as system of first order
with velocity variables $\Vo = \dot{\Po}$ and 
simultaneously enforces the constraints at velocity level (\ref{develo})
and the position constraints (\ref{constrigid}), where the latter are interpreted as
invariants and appended by means of extra Lagrange multipliers.

In this way, one obtains
an enlarged system
\begin{eqnarray}
            \dot{\Po} & = & \Vo - \Gr(\Po)^T \Lam \,, \nonumber \\
            \Mr(\Po) \, \dot{\Vo}& =& \fr(\Po,\Vo,t) - 
		      \Gr(\Po)^T \La  \,,  \label{ggl} \\
             \zero & = & \Gr(\Po) \, \Vo \,, \nonumber \\
             \zero & = & \gr(\Po) \, \nonumber
\end{eqnarray} 
with additional multipliers $\Lam(t) \in \mathbb{R}^{n_\Lac}$.
A straightforward calculation shows
\[ \zero = \frac{\D}{\D t} \, \gr(\Po) = \Gr(\Po)\dot{\Po}
      =  \Gr(\Po) \, \Vo -  \Gr(\Po)\Gr^T(\Po) \Lam = -  \Gr(\Po)\Gr^T(\Po) \Lam
\] 
and one concludes $\Lam = \zero$ since $\Gr(\Po)$ is of full rank and hence
$\Gr(\Po)\Gr^T(\Po)$ invertible.  With the 
additional multipliers $\Lam$ vanishing,
(\ref{ggl}) and the original equations of motion (\ref{f1}) coincide along any solution.
Yet, the index of the {\em GGL formulation\/}
(\ref{ggl}) is 2 instead of 3. 
Some authors refer to (\ref{ggl}) also as  
{\em stabilized index-2 system} \index{stabilized index-2 system}
 \cite{FuLe91}.

The last paragraphs on stabilized formulations of the Euler-Lagrange equations demonstrate that the development of theory and numerical methods for DAEs was strongly intertwined with the mathematical models. This holds for all application fields where DAEs arise. 
Even more, the application fields typically provide rich structural features of the model equations
that are crucial for determining the index and for the numerical treatment.

\section{Major results and numerical methods}\label{sect:dae-results}
Between 1989 and 1996, both theory and numerical analysis of DAEs were booming, and many groups in mathematics and engineering started to explore this new research topic. Driven by the development of powerful simulation packages in the engineering sciences, the demand for efficient and robust integration methods was growing steadily while at the same time, it had become apparent that higher index problems require stabilization measures or appropriate reformulations.

\subsection{Perturbation index and implicit Runge-Kutta methods}
The groundbreaking monograph on 
{\em The Numerical Solution of Differential-Algebraic Equations by {R}unge-{K}utta Methods} \cite{HaLR89} by Ernst Hairer, Christian Lubich and Michel Roche presented 
several new method classes, a new paradigm for the construction of convergence proofs, a new index concept, and the new RADAU5 code. 
From then on, Hairer and Lubich  played a very strong role in the further development of DAEs and corresponding numerical methods.

The {\em perturbation index\/} as defined in \cite{HaLR89} sheds a different light on DAEs and adopts the idea of a well-posed mathematical model.
While the differential index is based on successively differentiating
the original DAE (\ref{dae:fi}) until the obtained system can be solved for  $\dot{\xa}$,
the perturbation index 
measures the sensitivity of the solutions to perturbations in
the equation: \index{differential-algebraic equation!implicit}

The system $\FFa(\dot{\xa},\xa,t) = \zero$ has perturbation index $k \ge 1$ along a solution $\xa(t)$ on
$[t_0,t_1]$ if $k$ is the smallest integer such that, for all
functions $\hat{\xa}$ having a defect
\[ \FFa(\dot{\hat\xa},\hat\xa,t) = \pidy(t)\, , \]
there exists on
$[t_0,t_1]$  an estimate
\index{index!perturbation} 
\[ \| \hat{\xa}(t) - \xa(t) \| \le c 
\Big( \| \hat{\xa}(t_0) - \xa(t_0) \| + \max_{t_0 \le \xi \le t}
   \| \pidy(\xi) \| + \ldots + \max_{t_0 \le \xi \le t} \| \pidy^{(k-1)}(\xi) \| \Big) \]
whenever the expression on the right hand side is sufficiently small.
Note that the constant $c$ depends only on $\FFa$ and on the length of the interval, but not on the perturbation $\pidy$.
The perturbation index is $k=0$ if
\[
 \| \hat{\xa}(t) - \xa(t) \| \le c \Big( \| \hat{\xa}(t_0) - \xa(t_0) \| + \max_{t_0 \le \xi \le t}
    \| \int_{t_0}^{\xi} \pidy(\tau) \, \D\tau \| \Big)\,,
\]
which is satisfied for ordinary differential equations.

If the perturbation index exceeds $k=1$, derivatives
of the perturbation show up in the estimate and indicate
a certain degree of ill-posedness. E.g., if $\pidy$ contains
a small high frequency  term 
$\epsilon \sin \omega t$ with $\epsilon \ll 1$ and $\omega \gg 1$,
the resulting derivatives will induce a severe amplification in 
the bound for $\hat{\xa}(t) - \xa(t)$.

Unfortunately, the differential and the perturbation index are not equivalent
in general and may even differ substantially \cite{CG95}. 

The definition of the perturbation index is solely a prelude in \cite{HaLR89}. 
As the title says, most of the monograph deals with Runge-Kutta methods, in particular implicit ones. These are extended to linear-implicit systems $\Ea \dot{\xa} = \Fa(\xa,t)$ by discretizing $\dot{\xa} = \Ea^{-1} \Fa(\xa,t)$ and assuming for a moment that the matrix $\Ea$ is invertible. Multiplying the resulting 
scheme by $\Ea$, one gets the method definition
\index{Runge-Kutta methods}
\begin{subequations}\label{dae:irk}
\begin{eqnarray}
\Ea \Xa_i & = & \Ea \xa_0 + \hh \sum_{j=1}^s a_{ij} \Fa(\Xa_j,t_0+c_j \hh), \quad
i=1,\ldots,s; \\
\xa_1 & = & \left(1-\sum_{i,j=1}^s b_i \gamma_{ij}\right) \xa_0 + \hh \sum_{i,j=1}^s b_i \gamma_{ij}\Xa_j .
\end{eqnarray}
\end{subequations}
Here, the method coefficients 
are denoted by $(a_{ij})_{i,j=1}^s$ and $b_1, \ldots, b_s$ while $(\gamma_{ij})=(a_{ij})^{-1}$ is the inverse of the coefficient matrix,
with $s$ being the number of stages.
Obviously, (\ref{dae:irk}) makes sense also in the case where $\Ea$ is singular.

Using {\em stiffly accurate methods\/} for differential-algebraic equations is advantageous, which becomes evident if we consider the discretization of the semi-explicit system  
\begin{subequations}\label{dae:semi}
\begin{eqnarray}
      \dot{\ya} & = & \fa(\ya,\za), \label{dae:semia} \\
           \zero & = & \ga(\ya,\za) \label{dae:semib}
\end{eqnarray}
\end{subequations}
with differential variables $\ya$ and algebraic variables $\za$.
 The method (\ref{dae:irk}) 
then reads
\begin{subequations}\label{daes:irk}
\begin{eqnarray}
    \Ya_i  & = & \ya_0 + \hh \sum_{j=1}^s a_{ij} \fa(\Ya_j,\Za_j), \quad
i=1,\ldots,s,  \\
      \zero & =  & \ga(\Ya_i,\Za_i),
\end{eqnarray}
\end{subequations}
for the internal stages and 
\begin{subequations}\label{daes:irkup}
\begin{eqnarray}
\ya_1 & = & \ya_0 + \hh \sum_{j=1}^s b_j \fa(\Ya_j,\Za_j), \\
\za_1 & = & \left(1-\sum_{i,j=1}^s b_i \gamma_{ij}\right) \za_0 + \hh \sum_{i,j=1}^s b_i \gamma_{ij}\Za_j
\end{eqnarray}
\end{subequations}
as update for the numerical solution after one step. For stiffly accurate methods, we have $  \sum_{i,j=1}^s b_i \gamma_{ij} = 1$ and 
$\ya_1 = \Ya_s, \za_1 = \Za_s$. The update (\ref{daes:irkup}) is hence superfluous and furthermore, the constraint $\zero  =   \ga(\ya_1,\za_1)$
is satisfied by construction. 

It is not the purpose of this article to dive further into the world of Runge-Kutta methods, but like in numerical ODEs, the rivalry between multistep methods and Runge-Kutta methods also characterizes the situation for DAEs. While Linda Petzold's DASSL code is the most prominent multistep implementation, the RADAU5 and RADAU codes \cite{HaWa96,HaWa99}
represent the one-step counter parts and have also become widespread in various applications.

The competition for the best code was a major driving force in the numerical analysis of DAEs, and from time to time those in favor of multistep methods looked also at one-step methods, e.g., in   \cite{AsPe91}, and vice versa.

Simultaneously to the joint work with Ernst Hairer and Michel Roche,
Christian Lubich investigated a different class of discretization
schemes, the {\em half-explicit methods} \cite{Lubi89}.
These methods are tailored for semi-explicit DAEs and discretize the differential equations
explicitly while the constraint equations are enforced in an implicit 
fashion. 
As example, consider the Euler-Lagrange equations (\ref{f1}) with
velocity constraint (\ref{develo}). 
The half-explicit Euler method as generic
algorithm for the method class reads 
\begin{equation}\label{halbex} \begin{array}{rcl}
                 \Zs_{n+1} & = & \Zs_n + \hh \Vs_n \,,  \\
                 \MM(\Zs_n) \,\Vs_{n+1} & =& \MM(\Zs_n) \, \Vs_n + \hh \ff(\Zs_n,\Vs_n,t_n) - 
           \hh \GG(\Zs_n)^T \La_n \,, \\
                          \zero   & = & \GG(\Zs_{n+1}) \, \Vs_{n+1} \,.
\end{array} \end{equation}
Only a linear system of the form
\[
\left( \begin{array}{cc} \MM(\Zs_n) & \GG(\Zs_n)^T \\ \GG(\Zs_{n+1}) & \zero 
                        \end{array} \right) 
             \left( \begin{array}{c}  
                       \Vs_{n+1} \\ \hh \La_n  \end{array} \right) = 
             \left( \begin{array}{c} \MM(\Zs_n) \Vs_n + \hh \ff(\Zs_n,\Vs_n,t_n) \\
                      \zero \end{array} \right)
\]
arises here in each step. The scheme (\ref{halbex}) forms the basis for  
a class of extrapolation methods \cite{Lubi89,LNPE95}, and also for
 half-explicit Runge-Kutta methods as introduced in \cite{HaLR89}
and then further enhanced by Brasey \& Hairer \cite{BrHa93a} and
Arnold \& Murua \cite{ArMu98}.

These methods have in common that only information of the
velocity constraints is required. As remedy for the drift off, which
grows only linearly 
but might still be noticeable, 
the following projection, which is also due to Lubich \cite{Lubi89}, can be applied:
Let $\Zs_{n+1}$ and $\Vs_{n+1}$ denote the numerical solution of the system, obtained by integration from 
consistent values $\Zs_n$ and $\Vs_n$. 
Then, the projection consists of the following steps: 
\begin{subequations}\label{afw}
\begin{eqnarray}
\mbox{\quad}\hspace{-1.4cm}
& & \mbox{\em solve }    \left\{ \begin{array}{rcl}
             \zero & =& \MM(\tilde{\Zs}_{n+1}) (\tilde{\Zs}_{n+1}-\Zs_{n+1}) + \GG(\tilde{\Zs}_{n+1})^T \Lam, \\
             \zero & = & \gG(\tilde{\Zs}_{n+1})  \end{array} 
            \right. \, \mbox{\em for }
              \tilde{\Zs}_{n+1}, \Lam\,;  \label{afw1} \\[2mm]
\mbox{\quad}\hspace{-1.4cm}
& & \mbox{\em solve }    \left\{ \begin{array}{rcl}
             \zero & =& \MM(\tilde{\Zs}_{n+1}) (\tilde{\Vs}_{n+1}-\Vs_{n+1}) + \Gr(\tilde{\Zs}_{n+1})^T \ibv,   \\
       \zero & = & \Gr(\tilde{\Zs}_{n+1})\, \tilde{\Vs}_{n+1}  \end{array}
          \right. \, \mbox{\em for }
             \tilde{\Vs}_{n+1}, \ibv\,. \label{afw2} 
\end{eqnarray} 
\end{subequations}
A simplified Newton method can be used 
to solve  the nonlinear system (\ref{afw1}) while
(\ref{afw2}) represents a linear
system for $\tilde{\Vs}_{n+1}$ and $\ibv$ with similar structure.

The projection can also be employed for stabilizing the equations of motion 
with acceleration constraint (\ref{deaccel}) where the position and velocity constraints are invariants and not preserved by the time integration, see
Eich \cite{Eich93} and von Schwerin \cite{Schwe99}. Such
projection methods are particularly attractive in combination with explicit
ODE integrators. \index{projection methods}

\subsection{DAEs and differential geometry}

\begin{figure}[t]
\center
  \includegraphics[width=11cm]{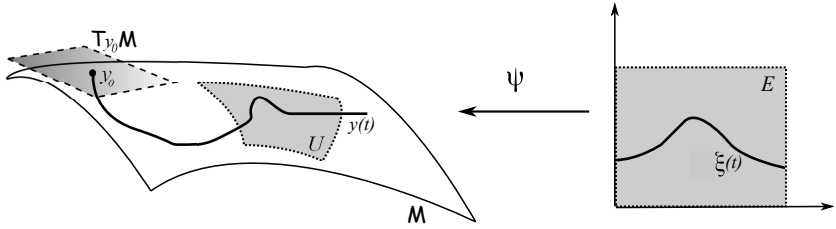}
   \caption{Manifold ${\cal M}$, tangent space ${\cal T}_\yac {\cal M}$,
    and local parametrization.}
   \label{locman}
\end{figure}

Already in 1984, Werner Rheinboldt had investigated DAEs from the viewpoint of differential geometry \cite{Rhei84}. While the approaches discussed so far
are mainly inspired by differential calculus and
algebraic considerations, a fundamentally different aspect comes into play by  his idea of {\em differential equations on manifolds}. 
\index{differential equation on manifold} 

Referring to \cite{AMR88,Arno80} for the theoretical underpinnings,
we shortly illustrate this approach by considering the semi-explicit system 
\index{differential-algebraic equation!semi-explicit} 
\begin{subequations}\label{dae:semi2}
\begin{eqnarray}
      \dot{\ya} & = & \fa(\ya,\za), \label{dae:semia2} \\
           \zero & = & \ga(\ya) \label{dae:semib2}
\end{eqnarray}
\end{subequations}
under the assumption
\begin{equation}\label{dae:vorind2}
\frac{\partial \ga}{\partial \ya}(\ya) \cdot
\frac{\partial \fa}{\partial \za}(\ya,\za) \in 
\mathbb{R}^{n_\zac \times n_\zac}
 \quad \mbox{is invertible}
\end{equation} 
in a neighborhood of the solution.
Clearly,  (\ref{dae:semi2}) is of index 2 where the constraint $\zero = \ga(\ya)$,
assuming sufficient differentiability,
defines the manifold \index{constraint!manifold}
\begin{equation}\label{dae:defmani}
  {\cal M} := \left\{ \ya \in \mathbb{R}^{n_\yac}: \ga(\ya) = \zero 
    \right\}  .
\end{equation}
The full rank condition (\ref{dae:vorind2}) for the matrix product
$\partial \ga / \partial \ya \cdot \partial \fa / \partial \za$ implies 
that the Jacobian $\Ga(\ya)=\partial \ga(\ya) / \partial \ya \in \mathbb{R}^{n_\zac \times
n_\yac}$ possesses also full rank $n_\zac$. Hence, for fixed 
$\ya \in {\cal M}$, the tangent space 
\index{tangent space}
\begin{equation}\label{dae:defmaniT}
 {\cal T}_\yac {\cal M} := \left\{ \yav \in \mathbb{R}^{n_\yac}: 
       \Ga(\ya) \yav = \zero 
    \right\}  
\end{equation}
is the kernel of $\Ga$ and has the
same dimension $n_\yac - n_\zac$ as the manifold ${\cal M}$.
Fig.~\ref{locman} depicts ${\cal M}$, ${\cal T}_\yac {\cal M}$,
and a solution of the DAE (\ref{dae:semi2}), which, starting 
from a consistent initial value, 
is required to proceed on the manifold.

The differential equation on the manifold ${\cal M}$ that is equivalent to 
the DAE  (\ref{dae:semi2}) is obtained as follows: 
The hidden constraint
\[   \zero = \Ga(\ya) \fa(\ya,\za)
\]
can be solved for $\za(\ya)$ according to the rank condition (\ref{dae:vorind2}) and the implicit function theorem.
Moreover, for $\ya \in {\cal M}$ it holds $\fa(\ya,\za(\ya)) \in
{\cal T}_\yac {\cal M}$, which defines a vector field on the manifold
${\cal M}$. Overall, 
\begin{equation}\label{dae:diffmani}
   \dot{\ya} = \fa(\ya,\za(\ya))  \quad \mbox{for } \ya \in {\cal M}
\end{equation}
represents then a differential equation on the manifold \cite{Rhei84}.

In theory, and also computationally \cite{manpak}, it is possible to
transform the differential equation (\ref{dae:diffmani}) from the manifold
to an ordinary differential equation in a linear space of 
dimension $n_\yac - n_\zac$. For this purpose, one introduces a 
{\em local parametrization} \index{local parametrization}
\begin{equation}\label{dae:deflocpara}
      \Par: \, E \rightarrow {\cal U}
\end{equation} 
where $E$ is an open subset of $\mathbb{R}^{n_\yac - n_\zac}$
and ${\cal U} \subset {\cal M}$, see Fig.~\ref{locman}.
Such a parametrization is not unique and holds only locally in general.
It is, however, possible to extend it to a family of parametrizations
such that the whole manifold is covered. For $\ya \in {\cal U}$ 
and local coordinates $\xlo \in E$ we thus get the relations
\[      \ya = \Par(\xlo), \qquad
        \dot{\ya} = \DPar(\xlo) \dot{\xlo}, \quad
        \DPar(\xlo) := \frac{\partial \Par}{\partial \xlo}(\xlo)
                        \in \mathbb{R}^{n_\yac \times (n_\yac - n_\zac)}.
\]
Premultiplying (\ref{dae:diffmani}) by the transpose of the Jacobian $\DPar(\xlo)$ 
of the parametrization and substituting $\ya$ by $\Par(\xlo)$, we 
arrive at 
\begin{equation}\label{dae:locssf}
    \DPar(\xlo)^T \DPar(\xlo) \dot{\xlo} = 
         \DPar(\xlo)^T \fa(\Par(\xlo),\za(\Par(\xlo))) \,.
\end{equation}
Since the Jacobian $\DPar$ has full rank for a 
valid parametrization, the matrix $\DPar^T \DPar$ is invertible, and
(\ref{dae:locssf}) constitutes the desired ordinary differential equation 
in the local coordinates $\xlo$. In analogy to a mechanical system 
in minimal coordinates, we call (\ref{dae:locssf})
a {\em local state space form}. 
\index{state space form} 

The process of transforming a differential equation
on a manifold to a local state space form constitutes
a {\em push forward\/} operator, while the reverse mapping
is called a {\em pull back\/} operator \cite{AMR88}.
It is important to realize that the previously defined
concept of an index does not appear in the theory of
differential equations on manifolds.
Finding hidden constraints by differentiation, however, 
is also crucial for the classification of DAEs from
a geometric point of view.

The geometrical viewpoint was also considered very early by Sebastian Reich 
\cite{reich1990}, but its full potential became clear only a couple of years later when the topic of geometric numerical integration emerged, cf.~\cite{HaLW02}.

\subsection{Singularly perturbed problems and regularization}

In the early days of DAEs, regularization was a quite popular means to convert the algebraic part into a differential equation. Motivated by physical examples such as stiff springs or parasitic effects in electric circuits, a number of authors have looked into this topic. 
Furthermore, it is also interesting to start with a singularly perturbed ODE, 
discretize it, and then to analyze the behavior of exact and numerical solution in the limit case.

To study an example for a semi-explicit system, we consider
 Van der Pol's equation 
\begin{equation}\label{dae:vdpo}
\epsilon \ddot \Poc + (\Poc^2-1) \dot \Poc + \Poc = 0 
\end{equation}
with parameter $\epsilon > 0$. This is an oscillator equation with a nonlinear damping term that acts as a controller. For large amplitudes $\Poc^2 > 1$, the damping term introduces dissipation into the system while for small values $\Poc^2 < 1$, the sign changes and the damping term is replaced by an excitation, leading thus to a self-exciting oscillator.  Introducing 
Li\'enhard's coordinates \cite{HaNW93}
\[   \zac := \Poc, \quad \yac := \epsilon \dot \zac 
                                     + (\zac^3/3 - \zac),
\] 
we transform (\ref{dae:vdpo}) into the first order system
\begin{subequations}\label{dae:vdp}
\begin{eqnarray} 
       \dot{\yac} & = & - \zac\, , \\
       \epsilon \dot{\zac} & = & \yac - \frac{\zac^3}{3} + \zac \,.
       \label{dae:vdpb}
\end{eqnarray} 
\end{subequations}
The case $\epsilon \ll 1$ is of special interest. In the limit
$ \epsilon = 0$, the equation  (\ref{dae:vdpb}) turns into a constraint and we arrive at the semi-explicit system
\begin{subequations}\label{dae:vdpr}
\begin{eqnarray} 
       \dot{\yac} & = & - \zac\, , \\
       0          & = & \yac - \frac{\zac^3}{3} + \zac \,. \label{dae:vdprb}
\end{eqnarray} 
\end{subequations}
In other words, Van der Pol's equation (\ref{dae:vdp}) in Li{\'e}nhard's coordinates is an example of a {\em singularly perturbed system\/}
\index{singularly perturbed system}
which tends to the semi-explicit \index{differential-algebraic equation!semi-explicit}  (\ref{dae:vdpr}) when
$\epsilon \rightarrow 0$.

Such a close relation between a singularly perturbed system and
a differen\-tial-algebraic equation is quite common and can be found 
in various application fields. Often, the parameter $\epsilon$ stands
for an almost negligible physical quantity or the presence of 
strongly different time scales. Analyzing the
{\em reduced system}, in this case (\ref{dae:vdpr}), usually proves successful to gain a better understanding of the original perturbed equation \cite{Mall74}. In the context of regularization methods,
this relation is also exploited, but in reverse
direction \cite{Hank90a}. One starts with a DAE such as (\ref{dae:vdpr}) and
replaces it by a singularly perturbed ODE, in this case (\ref{dae:vdp}).

In numerical analysis, the derivation and study of integration schemes via a singularly perturbed ODE has been termed the {\em indirect approach} \cite{HaLR89} and lead to much additional insight \cite{hairer1988error,LoPe86,Lubi93}, both for the differential-algebraic equation
as limit case and for the stiff ODE case. A particularly interesting method class for the indirect approach are Rosenbrock methods as investigated by Rentrop, Roche \& Steinebach \cite{ReRS89}.
\index{Rosenbrock methods} 

\subsection{General fully implicit DAEs}

In contrast to the solution theory in the linear constant coefficient case, the
treatment of fully implicit DAEs without a given internal structure is still challenging, even from today's perspective.
For this purpose, Campbell \cite{Ca93} introduced the derivative array as key
concept that carries all the information of the DAE system.
\index{derivative array}
The derivative array is constructed from the definition of the 
differential index, i.e., one  considers the equations
\begin{eqnarray}
		   \FFa(\dot{\xa},\xa,t) & = & \zero, \nonumber \\
		     \frac{ \D}{ \D t} \FFa(\dot{\xa},\xa,t)
		      = \frac{\partial}{\partial \dot{\xa}} \FFa(\dot{\xa},\xa,t)\, \xa^{(2)} + \ldots
		      & =& \zero, \label{dae:deriv} \\
				 & \vdots   &                 \nonumber      \\
		     \frac{ \D^k}{ \D t^k} \FFa(\dot{\xa},\xa,t)
		      = \frac{\partial}{\partial \dot{\xa}} \FFa(\dot{\xa},\xa,t)\, \xa^{(k+1)} + \ldots
		      & =& \zero \nonumber
\end{eqnarray}
for a DAE of index $k$. Upon discretization, (\ref{dae:deriv}) becomes an overdetermined system that can be tackled by least squares techniques.
The challenge in this procedure, however, is the in general unknown index $k$ and its determination.

Algorithms based on the derivative array are a powerful means for general unstructured DAE systems, and this holds even for the linear constant coefficient case since the computation of the Weierstrass form or the Drazin inverse are very sensitive to small perturbations and thus
problematic in finite precision arithmetic. For the derivative array, in contrast, so-called staircase algorithms have been developed that rely on orthogonal matrix multiplications and are much more stable \cite{BLMV2015}.

\subsection{Constrained Hamiltonian systems}\label{sect:hamilton}
\index{conservative system}
\index{Hamiltonian system}
In the conservative case, the Lagrange equations 
(\ref{f1}) of constrained mechanical motion
can be reformulated by the
transformation to Hamilton's canonical equations. 
This leads to a mathematical model that is typical for \emph{molecular dynamics}
simulations. Again, constraints come into play in this application field, and 
the time discretizations need to cope with the index and stability issues.

We define the Lagrange function \index{Lagrange!function}
\begin{equation}\label{defLTU}
 L(\Po,\dot{\Po}) := T(\Po,\dot{\Po}) - U(\Po)
\end{equation}
as the difference of kinetic and potential energy. \index{kinetic energy}
\index{potential energy}
The conjugate momenta $\Imp(t) \in \mathbb{R}^{n_\Poc}$
are then given by 
\begin{equation}\label{defmom}
   \Imp := \frac{\partial}{\partial \dot{\Po}}L(\Po,\dot{\Po}) = \Mr(\Po) \dot{\Po} \, , 
\end{equation}
and for the Hamiltonian we set
\begin{equation}\label{defHam}
   H := \Imp^T \dot{\Po} - L(\Po, \dot{\Po}) \,. 
\end{equation}
Since the velocity $\dot{\Po}$ can be expressed as $\dot{\Po}(\Imp,\Po)$ due to (\ref{defmom})
if the mass matrix is invertible,
we view the Hamiltonian as a function $H = H(\Imp,\Po)$.
Moreover, we observe that $H$ is the total energy of the system because
\[ H = \Imp^T \Mr(\Po)^{-1} \Imp - \frac{1}{2} \dot{\Po}^T \Mr(\Po) \dot{\Po} + U(\Po) = T + U \,.
\]
Using the least action principle (\ref{hamstarr})
in the new coordinates $\Imp$ and $\Po$ and applying
the Lagrange multiplier technique \index{Lagrange!multiplier} as above in the presence of constraints, 
we can express the 
equations of motion as 
\begin{eqnarray}
            \dot{\Po} & = & \frac{\partial}{\partial \Imp}H(\Imp,\Po)\,, \nonumber \\
            \dot{\Imp} & =& - \frac{\partial}{\partial \Po}H(\Imp,\Po) - 
	    		      \Gr(\Po)^T \La  \,, \label{hamcanon} \\ 
             \zero & = & \gr(\Po)  \,. \nonumber
\end{eqnarray} 
The {\em Hamiltonian equations\/} (\ref{hamcanon}) possess a rich mathematical structure 
that should be preserved by numerical methods. In the 1990s, this problem class led to
the new field of \emph{geometric integration}, see the monograph by Hairer, Lubich \& Wanner \cite{HaLW02} for an extensive exposition.
For brevity, we simply mention the \emph{SHAKE} scheme as one of the established
methods. In case of a Hamiltonian $H = \frac{1}{2}\Imp^T \Mr^{-1} \Imp  + U(\Po)$
with constant mass matrix, it reads
\begin{eqnarray}
    \Po_{n+1} - 2 \Po_n + \Po_{n-1} & = & -\tau^2 \Mr^{-1} (
                       \nabla U(\Po_n) + \Gr(\Po_n)^T \La_n ),  \label{shake} \\
            \zero & = & \gr(\Po_{n+1}).  \nonumber 
\end{eqnarray}
In each time thus a nonlinear system for $\Po_{n+1}$ and $\La_n$ needs to be solved,
which is closely related to the projection step (\ref{afw1}).

\section{Beyond classical DAEs}\label{sect:beyond}
At the end of the last century, new topics emerged that were
closely related to DAEs but also included important aspects from other fields. 
In particular, the topic of Partial Differential Algebraic Equations (PDAEs)
became attractive by then, driven by time-dependent partial differential equations
that were treated by the method of lines and that featured additional constraints.

\subsection{Navier-Stokes incompressible}
It requires convincing examples to demonstrate the benefits of a differential-algebraic 
viewpoint in the PDE context. One such example is sketched next.

A classical example for a PDAE is given by the Navier-Stokes equations
\begin{subequations}\label{nasto}
\begin{eqnarray}
        \dot{\svec{u}} + (\svec{u} \cdot \nabla)\svec{u} 
       + \mbox{$\frac{1}{\rho}$} \nabla p  & = & 
         \nu \svec{\Delta} \svec{u} + \svec{l},   \label{nastoa} \\
                0  & = & \nabla \cdot \svec{u} \label{nastob}
\end{eqnarray} 
\end{subequations}
for the velocity field $\svec{u}(x,t)$ and the pressure $p(x,t)$ 
in a $d$-dimensional domain $\Omega$, with 
mass density $\rho$, viscosity $\nu$, and source term $\svec{l}(x,t)$.
The second equation (\ref{nastob})
models the incompressibility of the fluid and defines a constraint for the velocity field.
For simplification, the convection term $ (\svec{u} \cdot \nabla)\svec{u}$ 
in (\ref{nastoa}) can be omitted, which makes the overall problem linear and more amenable 
for the analysis. 
In an abstract notation, the resulting Stokes problem then reads 
\begin{subequations} \label{nastoop}
\begin{eqnarray}
          \dot{\svec{u}} + \oA \svec{u} + \oB' p & = & \svec{l}, \\
                 \oB \svec{u} & = & 0 ,
\end{eqnarray} 
\end{subequations}
 with differential operators $\oA$ and $\oB$ expressing the Laplacian and the divergence, respectively. The notation $\oB'$ stands for the conjugate operator of $\oB$, which here is the gradient.

The discretization, e.g., by a
Galerkin-projection 
\[ \svec{u}(x,t) \doteq \svec{N}(x)\svec{q}(t), \qquad
                    p(x,t) \doteq \svec{Q}(x)\svec{\lambda}(t)
\]
with ansatz functions $\svec{N}$ and $\svec{Q}$ in some finite element spaces,
transforms the infinite-dimensional PDAE (\ref{nastoop}) to the DAE
\begin{subequations} \label{saddlep}
\begin{eqnarray}
          \svec{M} \dot{\svec{q}} + \svec{A} \svec{q} + \svec{B}^T \svec{\lambda} & = &  \svec{l} ,
                                      \\
                \svec{B} \svec{q}  & = &  \svec{0} .
\end{eqnarray}
\end{subequations}
While the mass matrix $\svec{M}$ and stiffness matrix $\svec{A}$ are symmetric positive definite and symmetric positive semi-definite, respectively, and easy to handle, the constraint matrix $\svec{B}$ is generated by mixing the discretizations for the velocity field and the pressure. It is well-known in mixed finite elements~\cite{Brezzi91} that a bad choice for the discretization will either result in a rank-deficient matrix $\svec{B}$ or in a situation where
the smallest singular value of $\svec{B}$ is approaching zero for a decreasing mesh size.
This means that the DAE~(\ref{saddlep}) may become singular or almost singular due to the spatial discretization.
The famous LBB-condition by Ladyshenskaja, Bab\v{u}ska, and Brezzi \cite{Brezzi91}
gives a means to classify the discretization pairs for $\svec{u}$ and $p$. If the matrix 
$\svec{B}$ has full rank, the index of the DAE (\ref{saddlep}) is $k=2$.

To summarize, PDEs with constraints such as the Navier-Stokes equations often feature a rich structure that should be exploited, and building on the available PDE methodology reveals interesting cross-connections with the differential-algebraic viewpoint. In this context,
the abstract formulation (\ref{nastoop}) as {\em transient saddle point problem} defines \index{transient saddle point problem}
a rather broad problem class where many application fields can be subsumed
\cite{simeon2013computational}.
 
By combining the state-of-the-art in DAEs with advanced PDE methodology and numerics,
powerful algorithms can then be developed that break new ground. Time-space adaptivitiy for PDAEs is one such topic where many different aspects are put together in order to set up numerical schemes with sophisticated error control. The work by Lang 
\cite{lang2013adaptive} defines a cornerstone in this field.
In  electrical circuit simulation, the inclusion of heating effects or semi-conductors results also in PDAE models where ODEs, DAEs, and PDEs are coupled via network approaches, see, e.g., \cite{Gu2001,ali2003elliptic}.

\subsection{Stochastic DAEs}
In the spring of 2006, Oberwolfach offered again the showcase for the latest developments in DAEs -- 25 years after the workshop where Bill Gear had first investigated the
mathematical pendulum (\ref{mathpend}) in Cartesian coordinates. The organizers were
Stephen Campbell, Roswitha M\"arz, Linda Petzold, and Peter Rentrop.
 Among the participants from all over the world was Bill Gear himself, and during the week it became evident that DAEs were now well-established in many fields.

The same year, the book by Kunkel \& Mehrmann \cite{KuMe06} appeared, which shed new light on topics such as boundary value problems in differential-algebraic equations and the numerical treatment of fully implicit systems (\ref{dae:fi}). 

Most talks at the meeting addressed the field of PDAEs, but
among the other prominent topics were also
optimization and optimal control problems with constraints described by DAEs, see, e.g.,
\cite{callies2008optimal,korkel2004numerical}
and model order reduction for descriptor systems \cite{reis2007stability}.
An emerging topic at the time were \emph{stochastic differential-algebraic equations} or SDAEs in short.
Since many models in science and engineering contain uncertain quantities, it is natural to extend the methodology for DAEs by corresponding random terms. This could either be a
parameter or coefficient that is only known approximately or even an extra diffusion term in the differential equation that is expressed in terms of a Wiener process. For the constant coefficient system (\ref{dae:constc}), 
such a diffusion term leads to the linear SDAE
\index{differential-algebraic equation!stochastic}
\begin{equation}\label{SDAE}
  \Ea {\D}{\xa}(t) + \Ba \xa(t) = \ca(t) + \svec{C} {\D} \svec{W}(t)
\end{equation}
with a Wiener process $\svec{W}$ in $\mathbb{R}^{n_x}$ and a square matrix $\svec{C}$. 
For work in this field and applications in electrical circuit analysis we 
refer to \cite{higham2002strong,winkler2003stochastic}.

\section{Nonsmooth dynamical systems} \label{sect:nonsmoothdynsys}

\begin{figure}
 \centering
 \subcaptionbox{A simple electrical circuit with an ideal diode $D$.\label{fig:circuit}}
 {\begin{tikzpicture}[
    circuit ee IEC,
    x = 3cm, y = 2cm,
    every info/.style = {font = \scriptsize},
    set diode graphic = var diode IEC graphic,
    set make contact graphic = var make contact IEC graphic,
  ]
  \draw (0,0) to [diode]     node[above,yshift=0.2cm] {\scriptsize $I\geq 0, V \geq 0$} 
                             node[below,yshift=-0.2cm] {$D$} (1,0) 
              to [capacitor] node[right,xshift=0.2cm] {$C$}  (1,1) 
              to [resistor]  node[above,yshift=0.2cm] {$R$}  (0,1) 
              to [inductor]  node[left,xshift=-0.2cm] {$L$}  (0,0);
\end{tikzpicture}}
 \subcaptionbox{A simulation of granular material with perfect unilateral contacts and friction\label{fig:granular}}
 {\includegraphics[width=0.62\textwidth]{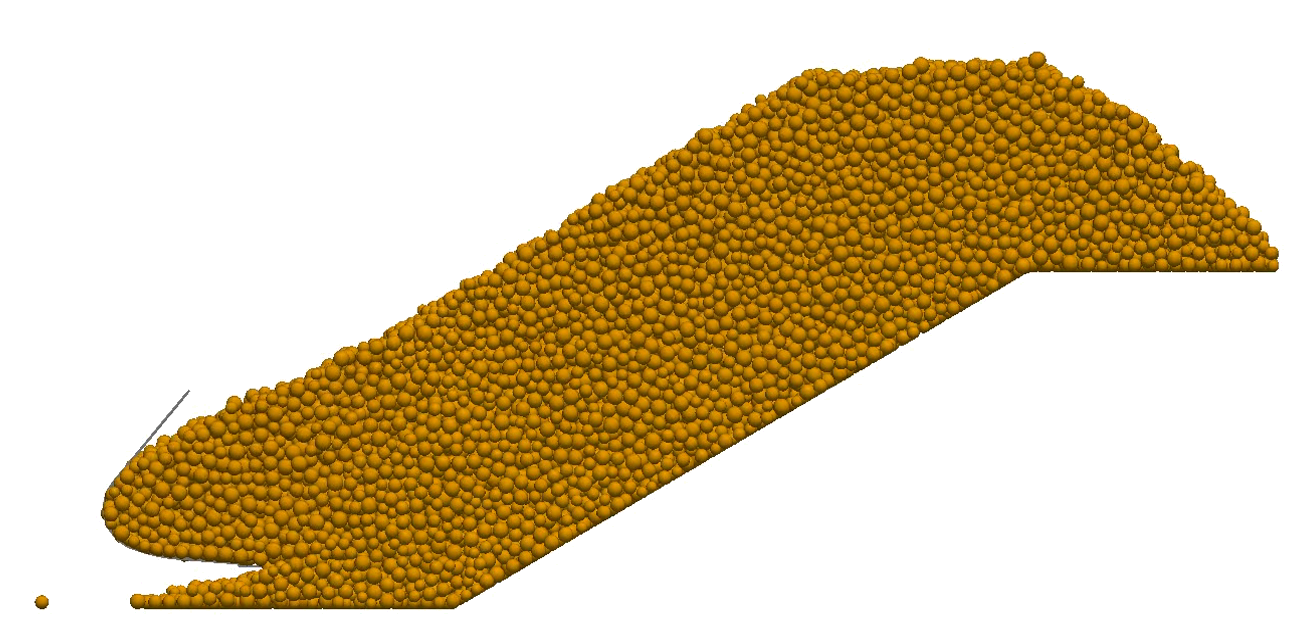}}
 \caption{Examples from electrical and mechanical engineering that include nonsmooth phenomena.}
\end{figure}
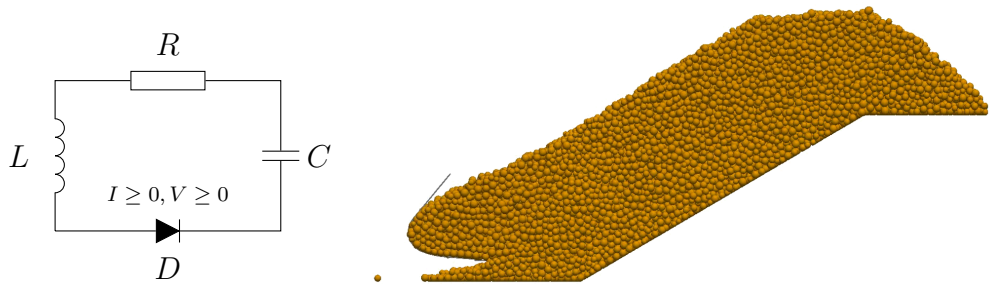

Differential equations for dynamical systems become particularly challenging when nonsmooth functions must be considered in one way or the other. Economical and financial mathematics have been a driving force for nonsmooth dynamical systems, where finding equilibria in market situations can be solved using variational inequalities that are closely related to nonlinear programming and convex optimization. Nonsmooth differential equations play a role in optimal control, robotic path planning and perturbation analysis. They occur also in electrical circuits with ideal diodes, see Figure~\ref{fig:circuit}. An ideal diode transmits electrical current only in one direction and blocks it in the other. If the electrical current $I$ is reversed, the conductivity of the ideal diode suddenly changes from 0 to $\infty$.  In other words,  the electric current at the diode must be positive at all times, $I(t) \geq 0$. If the current is nonzero, $I(t) > 0$ the voltage must be zero, $V(t)=0$. If there is a nonzero voltage $V(t) >0$ at the diode, this means, that the diode is blocking and does not transmit any current, $I(t) = 0$. This translates into the complementarity condition
\[ I(t) \geq 0, \hskip0.5cm V(t) \geq 0, \hskip0.5cm I(t)\cdot V(t) = 0. \]
at the diode. Which formalism can we use to incorporate inequalities and complementarity conditions into Kirchhoff's circuit laws? Another example are multibody dynamical systems with impacts, where contacts are modeled using a positivity constraint $g(\Po) \geq 0$ of a signed distance function, see for instance Figure~\ref{fig:granular}. Without long range and adhesive forces, the contact force $f_c$ between two contacting bodies must be nonnegative. In addition, it must be zero if the contact gap is positive. This again yields a complementarity problem
\[ g(\svec{q}) \geq 0, \hskip0.5cm f_c \geq 0, \hskip0.5cm g(\svec{q}) \cdot f_c = 0. \]

It is easy to construct examples, where the velocity of a rigid body must be discontinuous to enforce an inequality constraint. Think of a point mass accelerating under gravity towards the ground. At the moment of impact the velocity of the particle must jump from a negative to a nonnegative value instantaneously. This toy example demonstrates the two main difficulties that must be tackled in nonsmooth dynamical systems.
\begin{enumerate}
\item DAEs give us a mathematical structure to describe mechanical systems subject to \emph{algebraic} constraints. How can we incorporate \emph{inequality} constraints $g(\Po) \geq 0$?
\item Newton's second axiom 
\begin{equation}
\svec{M}\ddot{\qe} = \svec{f}(\qe,\dot{\qe},t) \label{eq:newton}
\end{equation}
implies that the trajectory $\qe(t)$ is at least once continuously differentiable, which is obviously not the case here. Clearly, Equation~\eqref{eq:newton} still has validity almost everywhere, but how can we handle the instances in time with impacts?
\end{enumerate}

In the following, we will take a closer look on the treatment of dynamical systems that are constrained to a feasible set:
\begin{align}
 &\FFa(\qe,\dot{\qe},t) = \boldsymbol{0}, \notag \\
 &\qe \in K = \left\{\; \qe \;\left| \; \svec{g}(\qe) \geq \boldsymbol{0} \; \right. \right\}
\end{align}
for some appropriate function $\svec{g}:\mathbb{R}^{n_q} \to \mathbb{R}^{n_{\lambda}}$. 

\subsection{A short zoology} \label{sect:dviHistory}

This section provides the most important concepts with regard to nonsmooth dynamical systems, including Moreau's sweeping process, differential inclusions, projected dynamical systems, variational inequalities, complementarity dynamic systems, differential variational inequalities and finally measure differential equations and measure differential inclusions. The list is far from complete and many more authors than the ones cited here have made important contributions. In addition, the overview is biased as the subject is approached from the perspective of nonsmooth mechanical problems.
For brevity, nonsmooth phenomena due to friction are excluded from this introduction all together.  

\subsubsection{Unilateral constraints and Moreau's sweeping process} \index{Moreau's sweeping process}

\begin{figure}[t]
 \centering
 \subcaptionbox{The tangent space $T_x\mathcal{M}$ and its annihilator $N_x\mathcal{M}$ at a point $x$ of a manifold $\mathcal{M}$, compare Figure~\ref{locman}.\label{fig:tangentspace}}
 {\begin{tikzpicture}[scale=0.66]

\pgfdeclarelayer{ft}
\pgfdeclarelayer{bg}
\pgfsetlayers{bg,main,ft}

\newcommand{\pgfextractangle}[3]{%
    \pgfmathanglebetweenpoints{\pgfpointanchor{#2}{center}}
                              {\pgfpointanchor{#3}{center}}
    \global\let#1\pgfmathresult  
}

\tikzfading 
[
  name=fade out,
  inner color=transparent!0,
  outer color=transparent!100
]

 \def\a{0.66};
 \def\b{0.4};
 \def\tmax{2};
 \def\N{50};
 
 \def\ta{-0.75};
 \def\tb{\tmax};

 \draw[color=black,line width=1.5,samples=\N,variable=\t,domain=-1.1*\tmax:1.1*\tmax]
    plot({\a*sinh(\t)},{-\b*cosh(\t)}) node[right] {$\mathcal{M}$};

 \def\r{1}
 \coordinate (O) at (0,0);
 \coordinate (pa) at ({\a*sinh(\ta)},{-\b*cosh(\ta)});
 \coordinate (ta) at ({\a*cosh(\ta)},{-\b*sinh(\ta)});
 \pgfextractangle{\angle}{O}{ta};
 
 \draw[line width=1.5, red] (pa) -- +(\angle-180:2*\r);
 \draw[line width=1.5, red] (pa) -- +(\angle:2*\r) node[right] {$T_{x}\mathcal{M}$};
 
 \draw[line width=1.5, blue] (pa) -- +(\angle-270:\r) node[above] {$N_{x}\mathcal{M}$};
 \draw[line width=1.5, blue] (pa) -- +(\angle-90:\r);
 
 \node[label={[xshift=-0.33cm,yshift=-0.4cm]\small $x$}] at (pa)  {$\bullet$};
 
\end{tikzpicture}} \hskip0.05\textwidth
 \subcaptionbox{The tangent cones $T_{x_i}K$ and normal cones $N_{x_i}K$ at three points $x_i$, $i=1,2,3$ in a set $K$. Note that $N_{x_3}K = \{\svec{0}\}$ and $T_{x_3}K$ takes up the whole space.\label{fig:tangentcone}}
 {\begin{tikzpicture}[scale=0.8]

\pgfdeclarelayer{ft}
\pgfdeclarelayer{bg}
\pgfsetlayers{bg,main,ft}

\newcommand{\pgfextractangle}[3]{%
    \pgfmathanglebetweenpoints{\pgfpointanchor{#2}{center}}
                              {\pgfpointanchor{#3}{center}}
    \global\let#1\pgfmathresult  
}

\tikzfading 
[
  name=fade out,
  inner color=transparent!0,
  outer color=transparent!100
]

 \def\a{0.66};
 \def\b{0.4};
 \def\c{0.66};
 \def\tmax{2};
 \def\N{50};
 
 \def\ta{-0.75};
 \def\tb{\tmax};

 
 \draw[dashed, color=black,line width=0.75,samples=\N,variable=\t,domain=-1.1*\tmax:1.1*\tmax]
    plot({\a*sinh(\t)},{-\b*cosh(\t)});
 \draw[dashed, color=black,line width=0.75,samples=\N,variable=\t,domain=-1.1*\tmax:1.1*\tmax]
    plot({\a*sinh(\t)},{\c*cosh(\t)-(\b+\c)*cosh(-\tmax)});
    
 \draw[name path = upper, color=black,line width=1.5,samples=\N,variable=\t,domain=-1.00*\tmax:1.00*\tmax]
    plot({\a*sinh(\t)},{-\b*cosh(\t)});
 \draw[name path = lower, color=black,line width=1.5,samples=\N,variable=\t,domain=-1.00*\tmax:1.00*\tmax]
    plot({\a*sinh(\t)},{\c*cosh(\t)-(\b+\c)*cosh(-\tmax)});
    
 \tikzfillbetween[of=lower and upper, on layer=bg] {gray!50};
 
 \node at (0,-1.8) {$K$};

 \def\r{1}
 \coordinate (O) at (0,0);
 \coordinate (pa) at ({\a*sinh(\ta)},{-\b*cosh(\ta)});
 \coordinate (ta) at ({\a*cosh(\ta)},{-\b*sinh(\ta)});
 \pgfextractangle{\angle}{O}{ta};
 
 \draw[line width=1.5, red] (pa) -- +(\angle-180:\r);
 \draw[line width=1.5, red] (pa) -- +(\angle:\r);
 \fill[red,path fading=south, opacity=0.5] (pa) ++(\angle-180:\r) arc ({\angle-180}:\angle:\r);
 
 \draw[line width=1.5, blue] (pa) -- +(\angle-270:\r) node[above] {$N_{x_1}K$};
 
 \node[label={[xshift=-0.33cm,yshift=-0.4cm]\small $x_1$}] at (pa)  {$\bullet$};
 
 \draw ($(pa)+(0.5,-0.5)$) -- +(1,1) node[above]  {$\textcolor{red}{T_{x_1}K}$};

 \def\r{1.5}
 \coordinate (pb) at ({\a*sinh(\tb)},{-\b*cosh(\tb)});
 \coordinate (tb1) at ({\a*cosh(\tb)},{-\b*sinh(\tb)});
 \coordinate (tb2) at ({\a*cosh(\tb)},{\c*sinh(\tb)});
 \pgfextractangle{\anglea}{O}{tb1};
 \pgfextractangle{\angleb}{O}{tb2};
 
 \draw[line width=1.5, red] (pb) -- +(\anglea+180:\r);
 \draw[line width=1.5, red] (pb)-- +(\angleb+180:\r);
 \fill[red, path fading=west, opacity=0.5] (pb) ++(\anglea+180:\r) arc ({180+\anglea}:{\angleb+540}:\r) -- (pb);

 \def\r{1};
 \draw[line width=1.5, blue] (pb) -- +(\anglea-270:\r);
 \draw[line width=1.5, blue] (pb) -- +(\angleb-90:\r);
 \fill[blue, path fading=east, opacity=0.5] (pb) ++(\anglea-270:\r) arc ({\anglea-270}:{\angleb-90}:\r) -- (pb);
 
 \node[label={[xshift=-0.1cm,yshift=0.00cm]\small $x_2$}] at (pb)  {$\bullet$};
 
 \node[label={[xshift=1.5cm,yshift=-0.5cm]$\textcolor{blue}{N_{x_2}K}$}] at (pb)  {};
 
 \draw ($(pb)+(-0.5,-0.25)$) -- +(1,-1) node[below]  {$\textcolor{red}{T_{x_2}K}$};
 
 \def\r{0.75}
 \coordinate (pc) at (-1,-2.25);

 \fill[red, path fading=fade out, opacity=0.5] (pc) circle(\r);

 \node[label={[xshift=-0.1cm,yshift=-0.25cm]\small $x_3$}] at (pc)  {$\textcolor{blue}{\bullet}$};
 
 \draw (pc) -- +(-0.25,-1) node[below] {$\textcolor{blue}{N_{x_3}K}$};

 \draw ($(pc)+(-0.25,0)$) -- +(-1,-0.33) node[left]  {$\textcolor{red}{T_{x_3}K}$};
 
\end{tikzpicture}}
 \caption{Comparison of tangent space and annihilator of manifolds to tangent and normal cones of sets.}
\label{fig:spaceVsCone}
\end{figure}
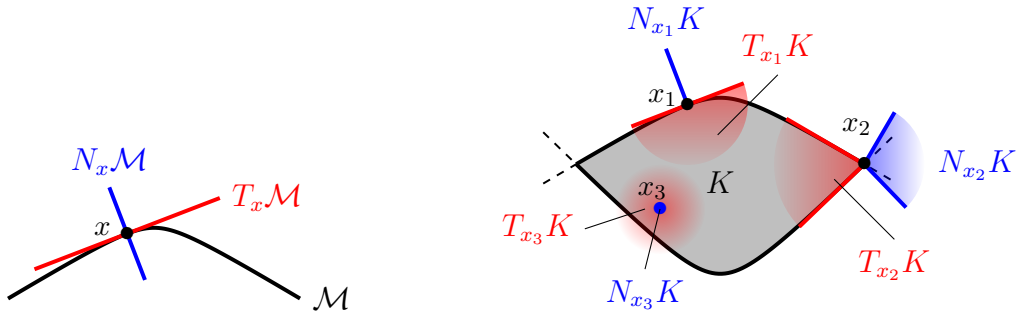

Jean Jacques Moreau studied mechanical systems subject to unilateral constraints since the mid 1960s \cite{moreau1966,moreau1977,moreau1988, moreau1990,moreau1993,moreau1999}. A \emph{unilateral constraint}\index{unilateral constraint} restricts the motion of a dynamical system only in one direction, e.g. by imposing inequality constraints. Imagine a particle's position $x(t)$ is restricted to a moving closed set $K(t)$. Consider a particular point in time $t^*$, where the  particle is located on the boundary of its feasible set, $x(t^*) \in \partial K(t^*)$. Assume for simplicity that the boundary can locally be described as a smooth manifold $\mathcal{M}$. Then we can find a tangent space $T_{x(t^*)}\mathcal{M}$ and a normal direction $\svec{n}$ to $T_{x(t^*)}\mathcal{M}$. The particle is only allowed to move freely in the half space to one side of the tangent space towards the interior of $K(t^*)$. In other words, the constraint restricts the motion \emph{unilaterally}. In contrast to this, if the particle were to be restricted to the manifold $\mathcal{M}$, the particle would not be allowed to move in either direction along the normal $\svec{n}$, i.e. it would be \emph{bilaterally constrained}. The space spanned by the normal vector is often called the annihilator of the tangent space. The assumption that the boundary of the moving set $K(t)$ can locally be described as a manifold is too restrictive in general. It suffices, if we can describe the boundary of the set through so--called tangent cones\footnote{See for instance \cite{rockafellar1998} for a definition of tangent cones of sets.\index{cone!tangent}} and normal cones rather than tangent spaces and their associated annihilators, see Figure~\ref{fig:spaceVsCone}.

A function $x$ is said to be a solution to the \emph{first order sweeping process} for the time--dependent set $K$, if
\begin{subequations} \label{eq:moreau}
\begin{align}
 x(0) &\in K(0), \\
 x(t) &\in K(t), \\
 -\dot{x}(t) &\in N_{x(t)}K(t), \label{eq:moreau-di}
\end{align}
\end{subequations}
where $N_{x(t)}K(t)$ denotes the normal cone to the set $K(t)$ at $x(t)$. It is defined as the polar of the tangent cone $T_{x(t)}K(t)$\index{cone!polar},
\[ N_{x(t)}K(t) = \left\{\left.\;  {\xi} \; \right| \; \langle {\xi}, {y} \rangle \leq 0 \hskip0.25cm \forall {y} \in T_{x(t)}K(t) \; \right\}. \] \index{cone!normal}%
The \emph{differential inclusion} (DI)\index{differential inclusion} in~\eqref{eq:moreau-di} means that if $x(t)$ is on the boundary of the set, the derivative must point inward, see Figure~\ref{fig:tangentcone}. In the interior of the set $K(t)$, it holds $N_{x(t)}K(t) = \{0\}$, and thus $\dot{x}(t) = 0$. In other words, the particle is \emph{swept} with the moving boundary.

\subsubsection{Projected dynamical systems} \index{projected dynamical system}

In \cite{nagurney1995}, a \emph{projected dynamical system} (PDS) is defined as
\begin{equation}
 \dot{x} = \Pi_K(x,-F(x)) \label{eq:pds}
\end{equation}
where
\[ \Pi_K(x,v) = \lim_{\delta \to 0} \frac{\proj_K(x+\delta v) -x}{\delta} \]
and $K$ is a closed convex set defined by constraints on the system. Here, $\proj_K(x)$ denotes the projection operator
\begin{equation}
 \proj_K(x)  = \arg\,\min\limits_{z \in K} \|x-z\|.\label{eq:proj}
\end{equation}
Because it holds
\[ \Pi_K(x,v) = \proj_{T_{x}K}(v), \]
Equation~\eqref{eq:pds} can be rewritten as a projection onto the tangent cone $T_xK$ \cite{brogliato2006}:
\begin{equation}
 \dot{x}=\proj_{T_{x}K}(-F(x)). \label{eq:pds2}
\end{equation}
Equations~\eqref{eq:pds} and~\eqref{eq:pds2} are first order ordinary differential equations with nonsmooth right hand sides. In Moreau's sweeping process, the direction of the derivative is prescribed only by the moving set $K$, while here, the direction of the derivative is mainly prescribed by the right hand side $-F(x)$. In contrast to Moreau's sweeping process, it holds $\dot{x} = -F(x) \neq 0$ in the interior of the set. The projection operator makes sure, that the system never moves towards the exterior of $K$ if $x$ is on the boundary, see Figure~\ref{fig:projectedsystem}. According to \cite{brogliato2006}, any solution to Equation~\eqref{eq:pds2} is also a solution of the differential inclusion
\begin{equation}
 -\dot{x}(t) \stackrel{\text{a.e.}}{\in} F(x(t)) + N_{x(t)}K, \label{eq:pds-di}
\end{equation}
which is Moreau's first order sweeping process, if $F(x) = 0$.

\begin{figure}[t]
 \centering
 \subcaptionbox{If $-F(x) \in T_xK$ it holds $\proj_{T_xK}(-F(x)) = -F(x)$ and therefore $\dot{x} = -F(x)$.\label{fig:projected-1}}
 {\begin{tikzpicture}[scale=0.8]

\pgfdeclarelayer{ft}
\pgfdeclarelayer{bg}
\pgfsetlayers{bg,main,ft}

\newcommand{\pgfextractangle}[3]{%
    \pgfmathanglebetweenpoints{\pgfpointanchor{#2}{center}}
                              {\pgfpointanchor{#3}{center}}
    \global\let#1\pgfmathresult  
}

\tikzfading 
[
  name=fade out,
  inner color=transparent!0,
  outer color=transparent!100
]

 \def\a{0.66};
 \def\b{0.4};
 \def\c{0.66};
 \def\tmax{2};
 \def\N{50};
 
 \def\ta{-0.75};
 \def\tb{\tmax};

 
 \draw[dashed, color=black,line width=0.75,samples=\N,variable=\t,domain=-1.1*\tmax:1.1*\tmax]
    plot({\a*sinh(\t)},{-\b*cosh(\t)});
 \draw[dashed, color=black,line width=0.75,samples=\N,variable=\t,domain=-1.1*\tmax:1.1*\tmax]
    plot({\a*sinh(\t)},{\c*cosh(\t)-(\b+\c)*cosh(-\tmax)});
    
 \draw[name path = upper, color=black,line width=1.5,samples=\N,variable=\t,domain=-1.00*\tmax:1.00*\tmax]
    plot({\a*sinh(\t)},{-\b*cosh(\t)});
 \draw[name path = lower, color=black,line width=1.5,samples=\N,variable=\t,domain=-1.00*\tmax:1.00*\tmax]
    plot({\a*sinh(\t)},{\c*cosh(\t)-(\b+\c)*cosh(-\tmax)});
    
 \tikzfillbetween[of=lower and upper, on layer=bg] {gray!50};
 
 \node at (-0.0,-1.8) {$K$};

 \def\r{1.5}
 \coordinate (O) at (0,0);
 \coordinate (pa) at ({\a*sinh(\ta)},{-\b*cosh(\ta)});
 \coordinate (ta) at ({\a*cosh(\ta)},{-\b*sinh(\ta)});
 \pgfextractangle{\angle}{O}{ta};
 
 \draw[line width=1.5, red] (pa) -- +(\angle-180:\r);
 \draw[line width=1.5, red] (pa) -- +(\angle:\r) node[right] {$T_{x}K$};
 \fill[red,path fading=south, opacity=0.5] (pa) ++(\angle-180:\r) arc ({\angle-180}:\angle:\r);

 \node[label={[xshift=-0.33cm,yshift=-0.4cm]\small $x$}] at (pa)  {$\bullet$};
 
 
 \draw[black!60!green,->, line width=2pt] (pa) -- +(1.5,-1) node[right] {$\dot{x}$};

\end{tikzpicture}} \hskip0.05\textwidth
 \subcaptionbox{If $-F(x) \notin T_xK$ it holds $\dot{x} = \proj_{T_xK}(-F(x)) \neq -F(x)$ and $-F(x) - \dot{x} \in N_xK$ \emph{(blue dotted line)}.\label{fig:projected-2}}
 {\begin{tikzpicture}[scale=0.8]

\pgfdeclarelayer{ft}
\pgfdeclarelayer{bg}
\pgfsetlayers{bg,main,ft}

\newcommand{\pgfextractangle}[3]{%
    \pgfmathanglebetweenpoints{\pgfpointanchor{#2}{center}}
                              {\pgfpointanchor{#3}{center}}
    \global\let#1\pgfmathresult  
}

\tikzfading 
[
  name=fade out,
  inner color=transparent!0,
  outer color=transparent!100
]

 \def\a{0.66};
 \def\b{0.4};
 \def\c{0.66};
 \def\tmax{2};
 \def\N{50};
 
 \def\ta{-0.75};
 \def\tb{\tmax};

 
 \draw[dashed, color=black,line width=0.75,samples=\N,variable=\t,domain=-1.1*\tmax:1.1*\tmax]
    plot({\a*sinh(\t)},{-\b*cosh(\t)});
 \draw[dashed, color=black,line width=0.75,samples=\N,variable=\t,domain=-1.1*\tmax:1.1*\tmax]
    plot({\a*sinh(\t)},{\c*cosh(\t)-(\b+\c)*cosh(-\tmax)});
    
 \draw[name path = upper, color=black,line width=1.5,samples=\N,variable=\t,domain=-1.00*\tmax:1.00*\tmax]
    plot({\a*sinh(\t)},{-\b*cosh(\t)});
 \draw[name path = lower, color=black,line width=1.5,samples=\N,variable=\t,domain=-1.00*\tmax:1.00*\tmax]
    plot({\a*sinh(\t)},{\c*cosh(\t)-(\b+\c)*cosh(-\tmax)});
    
 \tikzfillbetween[of=lower and upper, on layer=bg] {gray!50};
 
 \node at (-0.0,-1.8) {$K$};

 \def\r{2.5}
 \coordinate (pb) at ({\a*sinh(\tb)},{-\b*cosh(\tb)});
 \coordinate (tb1) at ({\a*cosh(\tb)},{-\b*sinh(\tb)});
 \coordinate (tb2) at ({\a*cosh(\tb)},{\c*sinh(\tb)});
 \pgfextractangle{\anglea}{O}{tb1};
 \pgfextractangle{\angleb}{O}{tb2};
 
 \draw[line width=1.5, red] (pb) -- +(\anglea+180:\r) node[midway,above,yshift=0.25cm] {$T_xK$};
 \draw[line width=1.5, red, name path = A] (pb)-- +(\angleb+180:\r);
 \fill[red, path fading=west, opacity=0.5] (pb) ++(\anglea+180:\r) arc ({180+\anglea}:{\angleb+540}:\r) -- (pb);

 
 \node at (pb)  {$\bullet$};
 \node[right,xshift=0.2cm] at (pb)  {$x$};
 
 \draw[black!60!green,->, line width=2pt] (pb) -- +(-0.33,-1.75) node[right,xshift=0.2cm] {$-F(x)$};
 \coordinate (F) at ($(pb)+(-0.33,-1.75)$);
 \path [name path = B] (F) -- +(\angleb+90:2);
 \path [name intersections={of=A and B,by=C}];
 \draw[line width=1.5, blue, dotted] (F) -- (C);
 \draw[black!60!green,->,line width=2] (pb) -- (C) node[midway,left,xshift=-0.25cm] {$\dot{x}$};

\end{tikzpicture}}
 \caption{Two example scenarios for projected dynamical systems.}
\label{fig:projectedsystem}
\end{figure}

Standard works on differential equations with discontinuous right hand side and differential inclusions are \cite{aubin1984,filippov1988}. Existence and uniqueness to such problems depends on properties of the set--valued map $K(t)$ and its tangent and normal cones. Nagurney and Zhang motivate their definition of projected dynamical systems as a means to unify the theory of dynamical systems and variational inequalities.

\subsubsection{Variational inequalities} \index{variational inequality}

Stuart Antman recites the emergence of variational inequalities in \cite{antman1983}. In 1959, Antonio Signorini posed the problem of finding the material displacements of a heavy deformable body resting on a rigid frictionless flat ground. This problem, now called the \emph{Signorini contact problem}, is particularly difficult, as the geometry of the contact region between the body and the ground is a priori unknown. Gaetano Fichera, a student of Signorini, studied existence and uniqueness of solutions to this problem using the calculus of variations and published his results in 1963 and 1964. In hindsight, Fichera's solution of Signorini's contact problem relies on the solution to a \emph{variational inequality} (VI). In a general setting, a variational inequality is defined as follows. Let $X$ be a Banach space, $K \subset X$ and $g:K \to X^*$ be a mapping from $K$ to the dual space $X^*$ of $X$. A variational inequality denotes the problem of finding $u \in K$, such that
\begin{equation}
 \langle g(u), u'-u \rangle \geq 0 \hskip0.5cm \forall\,u' \in K. \label{eq:vi}
\end{equation}
Figure~\ref{fig:vi} contains a geometric representation of a variational inequality in finite dimensions. Note that~\eqref{eq:vi} translates to $-g(u) \in N_uK$. The term variational inequality was first introduced by Stampacchia in 1965 and a theory around variational inequalities quickly grew during the second half of the 1960s. Most researchers on the subject at that time heavily cited Fichera's work. For more details on how Fichera's solution to Signorini's problem can be solved using variational inequalities and how his work influenced their analysis in the years to come, see \cite{antman1983}.

\subsubsection{Complementarity dynamic systems} \index{complementarity dynamic system}

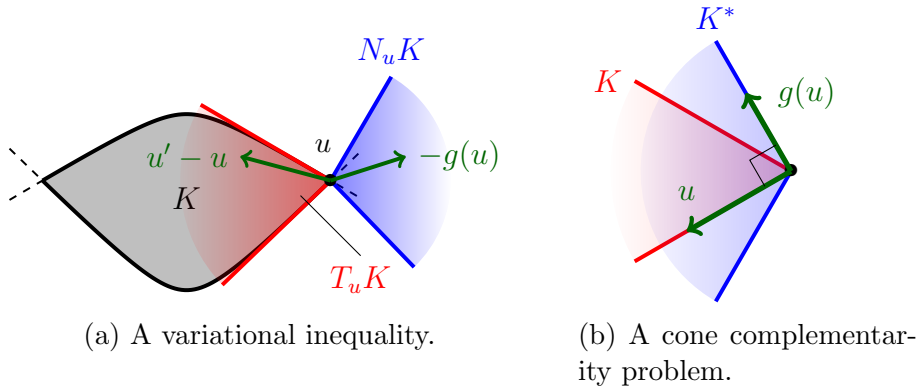
\begin{figure}[t]
 \centering
 \subcaptionbox{A variational inequality.\label{fig:vi}}
 {\begin{tikzpicture}[scale=0.8]

\pgfdeclarelayer{ft}
\pgfdeclarelayer{bg}
\pgfsetlayers{bg,main,ft}

\newcommand{\pgfextractangle}[3]{%
    \pgfmathanglebetweenpoints{\pgfpointanchor{#2}{center}}
                              {\pgfpointanchor{#3}{center}}
    \global\let#1\pgfmathresult  
}

\tikzfading 
[
  name=fade out,
  inner color=transparent!0,
  outer color=transparent!100
]

 \def\a{0.66};
 \def\b{0.4};
 \def\c{0.66};
 \def\tmax{2};
 \def\N{50};
 
 \def\ta{-0.75};
 \def\tb{\tmax};

 
 \draw[dashed, color=black,line width=0.75,samples=\N,variable=\t,domain=-1.1*\tmax:1.1*\tmax]
    plot({\a*sinh(\t)},{-\b*cosh(\t)});
 \draw[dashed, color=black,line width=0.75,samples=\N,variable=\t,domain=-1.1*\tmax:1.1*\tmax]
    plot({\a*sinh(\t)},{\c*cosh(\t)-(\b+\c)*cosh(-\tmax)});
    
 \draw[name path = upper, color=black,line width=1.5,samples=\N,variable=\t,domain=-1.00*\tmax:1.00*\tmax]
    plot({\a*sinh(\t)},{-\b*cosh(\t)});
 \draw[name path = lower, color=black,line width=1.5,samples=\N,variable=\t,domain=-1.00*\tmax:1.00*\tmax]
    plot({\a*sinh(\t)},{\c*cosh(\t)-(\b+\c)*cosh(-\tmax)});
    
 \tikzfillbetween[of=lower and upper, on layer=bg] {gray!50};
 
 \node at (0,-1.8) {$K$};

 \def\r{2.5}
 \coordinate (pb) at ({\a*sinh(\tb)},{-\b*cosh(\tb)});
 \coordinate (tb1) at ({\a*cosh(\tb)},{-\b*sinh(\tb)});
 \coordinate (tb2) at ({\a*cosh(\tb)},{\c*sinh(\tb)});
 \pgfextractangle{\anglea}{O}{tb1};
 \pgfextractangle{\angleb}{O}{tb2};
 
 \draw[line width=1.5, red] (pb) -- +(\anglea+180:\r);
 \draw[line width=1.5, red] (pb)-- +(\angleb+180:\r);
 \fill[red, path fading=west, opacity=0.5] (pb) ++(\anglea+180:\r) arc ({180+\anglea}:{\angleb+540}:\r) -- (pb);

 \def\r{2};
 \draw[line width=1.5, blue] (pb) -- +(\anglea-270:\r) node[above] {$N_uK$};
 \draw[line width=1.5, blue] (pb) -- +(\angleb-90:\r);
 \fill[blue, path fading=east, opacity=0.5] (pb) ++(\anglea-270:\r) arc ({\anglea-270}:{\angleb-90}:\r) -- (pb);
 
 \node[label={[xshift=-0.1cm,yshift=0.00cm]\small $u$}] at (pb)  {$\bullet$};
 
 \draw ($(pb)+(-0.5,-0.25)$) -- +(1,-1) node[below]  {$\textcolor{red}{T_{u}K}$};
 
 \draw[black!60!green,->,line width = 1.5pt] (pb) -- +(-1.5,0.4) node[left] {$u' - u$};
 \draw[black!60!green,->,line width = 1.5pt] (pb) -- +(1.25,0.4) node[right] {$-g(u)$};
 
\end{tikzpicture}} \hskip0.05\textwidth
 \subcaptionbox{A cone complementarity problem.\label{fig:ccp}}
 {\begin{tikzpicture}[scale=0.8]

\pgfdeclarelayer{ft}
\pgfdeclarelayer{bg}
\pgfsetlayers{bg,main,ft}

\newcommand{\pgfextractangle}[3]{%
    \pgfmathanglebetweenpoints{\pgfpointanchor{#2}{center}}
                              {\pgfpointanchor{#3}{center}}
    \global\let#1\pgfmathresult  
}

\tikzfading 
[
  name=fade out,
  inner color=transparent!0,
  outer color=transparent!100
]

 \def\a{0.66};
 \def\b{0.4};
 \def\c{0.66};
 \def\tmax{2};
 \def\N{50};
 
 \def\ta{-0.75};
 \def\tb{\tmax};

 
%
%
%

 \def\r{3} 
 \draw[line width=1.5, red] (pb) -- +(150:\r) node[left] {$K$};
 \draw[line width=1.5, red] (pb)-- +(210:\r);
 \fill[red, path fading=west, opacity=0.25] (pb) ++(150:\r) arc (150:210:\r) -- (pb);

 \def\r{2.5};
 \draw[line width=1.5, blue] (pb) -- +(120:\r) node[above] {$K^*$};
 \draw[line width=1.5, blue] (pb) -- +(240:\r);
 \fill[blue, path fading=west, opacity=0.25] (pb) ++(120:\r) arc (120:240:\r) -- (pb);
 
 \node at (pb)  {$\bullet$};
 
 \draw[black!60!green,->,line width = 2pt] (pb) -- +(210:2) coordinate (u) node[above,yshift=0.25cm] {$u$};
 \draw[black!60!green,->,line width = 2pt] (pb) -- +(120:1.5) coordinate (gu) node[right,xshift=0.25cm] {$g(u)$};
 
 \draw ($(pb)!0.5cm!(u)$) coordinate (t) -- ($(t)!0.5cm!-90:(u)$) -- ($(pb)!0.5cm!(gu)$);
 
\end{tikzpicture}}
 \caption{Geometric interpretations of a variational inequality and a cone complementarity problem in finite dimensions.}
\label{fig:vi-ccp}
\end{figure}

A convex cone $K$ is a set that is closed with respect to positive linear combinations, i.e. if $x$ and $y$ are in the cone $K$ and $\alpha, \beta \geq 0$, it must hold that $\alpha x + \beta y \in K$\index{cone!convex}. The variational inequality~\eqref{eq:vi} is equivalent to a \emph{cone complementarity problem} (CCP) \index{cone complementarity problem}
\begin{equation}
 g(u) \in K^*, \hskip0.5cm u \in K, \hskip0.5cm \langle g(u),u \rangle = 0,\label{eq:ccp}
\end{equation}
if the set $K$ is a convex cone, which is proved in \cite{karamardian1971}. Here,
\[ K^* = \left\{\; y \in X^* \;\left|\; \langle y, x \rangle \geq 0 \; \forall x \in K \; \right. \right\} \subset X^* \]
is the \emph{dual cone} to $K$\index{cone!dual}. Notice that the dual cone is the negative normal cone to $K$ at $0$. Because of the complementarity condition in~\eqref{eq:ccp} it is clear, that if either $g(u)$ or $u$ is in the interior of their respective cones $K^*$ and $K$, the other of the two must be zero (see Figure~\ref{fig:ccp}), hence the term \emph{complementarity}. 

The positive orthant $\mathbb{R}_+^n$ is such a convex cone and it holds $(\mathbb{R}_+^n)^* = \mathbb{R}_+^n$, i.e. it is self--dual. For this specific cone the CCP~\eqref{eq:ccp} resembles in part the well--known Karush--Kuhn--Tucker necessary conditions for optimization problems with inequality constraints $\svec{g}(\svec{u}) \geq 0$:
\[ \svec{g}(\svec{u}) \geq \svec{0}, \hskip0.5cm \svec{u} \geq \svec{0}, \hskip0.5cm \svec{u}^T\svec{g}(\svec{u}) = 0, \]
where the inequality ``$\geq$'' is to be understood componentwise.

\emph{Complementarity dynamic systems} are differential equations coupled to a complementarity problem. They were considered by several authors, including J. J. Moreau himself, to tackle mechanical systems with inequality constraints \cite{moreau1966,Loe82,Anitescu1997a,Jean1999,pfeiffer2003,Pfeiffer2006}. Their works were heavily influenced by results from nonlinear programming and optimization.

\subsubsection{Differential variational inequalities} \index{variational inequality!differential}

In \cite{aubin1984} a special version of a differential inclusion is mentioned, called a \emph{differential variational inequality}. Given a closed convex subset $K \subset X$ of a vector space $X$ and a set--valued map $F: K \to X$, a differential variational inequality consists of finding a function $x: [0,T] \to X$ satisfying
\begin{subequations}\label{eq:aubin-dvi}
\begin{align}
 x(t) &\in K  &\forall\;t \in [0,T], \\
 \dot{x}(t) &\in F(x(t)) - N_{x(t)}K &\text{for a.a.}\;t\in [0,T]. \label{eq:aubin-dvi-di}
\end{align}
\end{subequations}
Notice that the differential inclusion~\eqref{eq:aubin-dvi-di} coincides with~\eqref{eq:pds-di}. The relation to variational inequalities becomes clear, if we write $\dot{x}(t) = f(x(t)) + g(t)$ for $f(x(t)) \in F(x(t))$ and $g(t) \in -N_{x(t)}K$. With this, \eqref{eq:aubin-dvi} reads
\begin{subequations}\label{eq:aubin-dvi-2}
\begin{align}
 &x(t) \in K  &\forall\;t \in [0,T], \\
 &\langle \dot{x}(t) - f(x(t)), y - x(t) \rangle \geq 0 &\text{for a.a.}\;t\in [0,T] \text{ and all } y \in K.
\end{align}
\end{subequations}

A different, slightly more general, definition of a differential variational inequality comes from Pang and Stewart \cite{pang2008}. They define a finite dimensional differential variational inequality as a differential equation coupled to a variational inequality
\begin{subequations} \label{eq:pang-dvi}
\begin{align}
 &\dot{\svec{x}} = \svec{f}(t,\svec{x}(t),\svec{u}(t)), \\
 &\langle \svec{F}(t,\svec{x}(t),\svec{u}(t)), \svec{u}' - \svec{u} \rangle \geq 0 \hskip0.5cm \forall \svec{u}' \in K, \label{eq:pang-dvi-vi} \\
 &\Gamma(\svec{x}(0),\svec{x}(T)) = 0. \label{eq:pang-dvi-bc}
\end{align}
\end{subequations}
The problem consists of finding functions $\svec{x}:[0,T] \to \mathbb{R}^n$ and $\svec{u}:[0,T] \to K \subset \mathbb{R}^m$ satisfying Equation~\eqref{eq:pang-dvi} given the continuous functions $\svec{f}$, $\svec{F}$, $\Gamma$ and a subset $K \subset \mathbb{R}^m$. In this finite dimensional setting the dual pairing~\eqref{eq:pang-dvi-vi} is the Euclidean scalar product and Equation~\eqref{eq:pang-dvi-bc} is a prescribed initial or boundary condition. If $K$ is a convex cone, the variational inequality~\eqref{eq:pang-dvi} can be expressed as a cone complementarity problem. To differentiate their definition~\eqref{eq:pang-dvi} from the definition~\eqref{eq:aubin-dvi} given by Aubin and Cellina, Pang and Stewart call problems of the type~\eqref{eq:aubin-dvi} \emph{variational inequalities of evolution}  \index{variational inequality!of evolution}(VIE), a convention that is adapted here. Pang and Stewart show in the same article, that the VIE is a special case of a DVI and suggest a unified version of the two, that essentially consists of a differential equation, an algebraic equation and a variational inequality. Finally, it should be noted that if the set $K$ is a convex cone, a DVI and a VIE are conceptually equivalent. 

Let $\mathcal{U}$ denote the set of solutions to the variational inequality~\eqref{eq:pang-dvi-vi}. Then, the DVI~\eqref{eq:pang-dvi} can be rewritten as differential inclusion
\[ \dot{\svec{x}} \in \svec{f}(t,\svec{x}(t),\mathcal{U}), \]
The DVI~\eqref{eq:pang-dvi} can also be rewritten as a  DAE: The function $\svec{u}$ solves the variational inequality~\eqref{eq:pang-dvi-vi} if and only if
\[ \svec{0} = \svec{g}(t,\svec{x},\svec{u})=\svec{u} - \proj_K(\svec{u}-\svec{F}(t,\svec{x},\svec{u})),\]
where $\proj_K$ is the orthogonal projection onto the set $K$ as defined in \eqref{eq:proj}. Thus, the variational inequality can be replaced by a nonsmooth algebraic constraint. Pang and Stewart plea their case why DVIs deserve a special treatment, even if they can be considered as either DIs or DAEs. They argue, that the DI theory is too general to be of practical use, since many of its assumptions are cumbersome to prove for specific situations. While the theory of DAEs appeals to them, it is problematic as the differentiability of the algebraic constraints play an important role in it. If the algebraic constraint however is nonsmooth, as is the case if it is constructed from a variational inequality, many results cannot be used.

Still, it is obvious that the authors of~\cite{pang2008} were inclined to adopt as much terminology from and analogy to DAEs for their analysis as possible. As for DAEs, the unknowns of a DVI are separated into differential variables $\svec{x}$ and algebraic variables $\svec{u}$. While a DAE can be seen as an ODE coupled to an algebraic equation through the algebraic variables, a DVI is a differential equation coupled to a variational inequality. Even the concept of the index of a DVI is adopted. A DVI of index zero is just an ordinary differential equation. Under certain conditions for the set--valued map $\svec{F}$, the algebraic unknown $\svec{u}$ can be written as a function of $\svec{x}$ and $t$ using the implicit function theorem. One differentiation then reveals the DVI as system of differential equations in $\svec{x}$ and $\svec{u}$. This is called a DVI of index one. According to the authors, higher index DVIs can be considered as well, but they restrict their analysis to DVIs of index one with absolutely continuous solutions.

\subsubsection{Derivatives of functions of bounded variation}

Differential variational inequalities and the related concepts help us to describe dynamical systems subject to inequality constraints. We opened the section on nonsmooth dynamical systems with another problem. If the velocity of a mechanical system subject to unilateral constraints is discontinuous, how do we interpret accelerations and forces? The results in the next two sections are taken from \cite{natanson1977,moreau88b,glocker2000,acary2008,gerdts2012optimal}.

For absolutely continuous functions $x:[0,T] \to \mathbb{R}$, the derivative $v = \dot{x}$ is of bounded variation. This entails that $v$ is smooth almost everywhere, but can have a countable number of finite jump discontinuities, at which $v=\dot{x}$ is not defined per se\footnote{We could easily define the function $v(t_c)$ at the discontinuity points $t_c$ as either the left limit $v^-(t_c) = \lim\limits_{h\to 0}v(t_c - h)$ or the right limit $v^+(t_c) = \lim\limits_{h\to 0}v(t_c + h)$, with $h > 0$, as convention. While this potentially helps the analysis, it does not change the structure of the problem.}. Still, we can write
\[ x(t) = x(0) + \int_0^t \dot{x}(\tau) \;\D \tau, \]
and just jump over the discontinuity points in the interval $[0,t]$, as they are a zero--set of the Lebesgue measure $\D \tau$ anyway.

How do we define the derivative of a function of bounded variation $v$? It turns out, that measure theory helps. If the function can have jumps, it changes infinitely fast from its left limit to its right limit at such a discontinuity point. The discontinuity points comprise only a zero--set of the Lebesgue measure $\D t$, but surely the infinitely fast changes in $v$ cannot simply be ignored. We must capture the changes of the function using measures, for which the set of discontinuity points is not a zero set. 

We can construct the so--called Lebesgue--Stieltjes measure\index{Lebesgue--Stieltjes measure} $\D v$ from a function of bounded variation $v$, so that it holds
\[ v(t) = v(0) + \int_0^t \mathrm{d}v. \]
In this sense, the Lebesgue--Stieltjes measure $\mathrm{d}v$ plays the role of the derivative of $v$. The measure captures the changes in $v$ regardless of whether they happen smoothly or abruptly. Using the Lebesgue decomposition and Radon--Nikodym theorems, the Lebesgue--Stieltjes measures of $v$ can be decomposed into three parts,
\[ \D v = a(t) \D t + \D \delta_p + \D \mu_s. \]
Here, $a(t)$ is the so--called Radon--Nikodym density\index{density!Radon--Nikodym} of the continuous part of the measure $\D v$ with respect to the Lebesgue measure $\D t$, $\D \delta_p$ is a sum of Dirac delta measures that captures all jumps of $v$ and $\D \mu_s$ is a singular measure associated with the Cantor part of $v$. If $v$ has no Cantor part, it is called a function of specially bounded variation and the last term vanishes. This is a valid assumption for a large class of problems. If in addition it is absolutely continuous, the second term vanishes and it holds $\D v = a(t) \D t$ with $a = \dot{v}$. 

Note that the Lebesgue--Stieltjes measure $\D v$ corresponds to the distributional derivative of $v$:
\[ \langle \mathrm{D} v, \varphi \rangle = -\int_0^T v \dot{\varphi}\;\D t = -\int_0^T v \;\D\varphi = \int_0^T \varphi \;\D v \hskip0.5cm \forall \varphi \in C_{0}^\infty(0,T), \]
because all functions in $C_{0}^\infty(0,T)$ are absolutely continuous and thus have a density w.r.t. $\D t$. Here we use integration by parts for Lebesgue--Stieltjes integrals:
\begin{equation}
 \int_0^T v\;\D w = \left[vw \right]_0^T - \int_0^T w\;\D v \label{eq:intbyparts}
\end{equation}
for all $v,w$ with bounded variation.

\subsubsection{Measure differential equations and measure differential inclusions}

A \emph{measure differential equation} (MDE) \index{measure differential equation} takes the form
\begin{equation}
 \D x = F(t,x)\;\D t + G(t,x) \;\D u, \label{eq:mde}
\end{equation}
where $\D x$ and $\D u$ denote the Lebesgue--Stieltjes measures of the functions of bounded variation $x$ and $u$  and $\D t$ denotes the Lebesgue measure of time \cite{leela1974, pandit1982}. The corresponding initial value problem consists of finding $x:[0,T] \to \mathbb{R}$ given an initial value $x(0)=x_0$. Measure differential equations were motivated from Optimal Control Theory and the theory of perturbed systems. Here, the function $u$ takes the role of inputs or external perturbations. In control theory, $u$ drives the state $x$ of the system to a certain target. In certain applications it might be useful to allow jumps in the input, which in turn result in jumps in the state and vice-versa. In perturbation theory, impulsive perturbations are studied. See \cite{piccoli2017} for a recent, slightly different approach towards measure differential equations in the context of differential geometry.

One way to read Equation~\eqref{eq:mde} is as a short-form for the variational condition
\begin{equation}
 \int_0^T \varphi \;\D x = \int_0^T \varphi \;F(t,x) \;\D t + \int_0^T \varphi\; G(t,x)\;\D u \hskip0.5cm \forall \varphi \in \mathcal{T}, \label{eq:mde-variational}
\end{equation}
for some appropriate space of test functions $\mathcal{T}$. If $u$ were absolutely continuous, we could write $\D u = \dot{u} \;\D t$. Then it is obvious from~\eqref{eq:mde}, that $\D x$ has a density $F(t,x) + G(t,x)\dot{u}$ with respect to $\D t$ and thus we can write $\D x = \dot{x}\;\D t$. Finally, using the fundamental lemma of calculus of variations, we could derive a differential equation
\[ \dot{x} = F(t,x) + G(t,x) \dot{u} \]
from~\eqref{eq:mde-variational}.

Moreau was well aware that the sweeping process results in discontinuous velocities. He introduced \emph{measure differential inclusions} (MDI)\index{differential inclusion!differential} as an extension to the differential inclusion in~\eqref{eq:moreau} \cite{moreau1988, monteiromarques1993,glocker2000,2001stewart,paoli2005}. In \cite{monteiromarques1993} a \emph{measure differential} inclusion is defined as
\[ \D u \in N_{u}K, \]
for a convex subset $K \subset X$ of a Hilbert space $X$. It should be noted here, that this definition would have just as much meaning in a Banach space, which will be exploited in the next section. In this abstract setting, the normal cone is a subset of the dual space $X^*$ of $X$. In the same book, existence and uniqueness results are given for second order measure differential inclusions that result from mechanical systems with unilateral constraints and impacts. See also \cite{acary2008b} for higher order MDIs in the context of Moreau's sweeping process.

Up to this point all considerations concerning normal and tangent cones could be directly translated to the finite dimensional spaces in which the functions take their values. In other words, if we consider functions taking values in $\mathbb{R}^n$, the dual pairings $\langle \cdot,\cdot \rangle$ can be interpreted as the Euclidean scalar product in $\mathbb{R}^n$ and the negativity condition in the definition of the normal cone must hold in a pointwise sense. This gives us a geometric interpretation of tangent and normal cones, complementarity problems and variational inequalities. As soon as distributional derivatives play a role, this interpretation becomes more and more difficult and the problems are usually considered in infinite dimensional function spaces.

\subsubsection{A very diverse field}

In the past 60 years, many researchers devoted themselves to dynamical systems that encounter nonsmooth phenomena. These researchers borrowed results from convex and functional analysis, numerical analysis, nonlinear programming and optimization, measure theory and differential geometry to produce a rich structure and a diverse selection of concepts to describe ``\emph{nonsmoothness}'' of differential equations. DAEs did not play a major role in the early development of the field. It is only recently that the DAE theory has become an important influence on the field of nonsmooth dynamical systems.

Many of introduced concepts are very similar and under mild assumptions some can be shown to be equivalent \cite{brogliato2006}. The diversity of the field comes at a small price however, which is the lack of a common language. It can be very hard to judge which of these concepts is most appropriate for the application at hand. It is extremely difficult to find the vocabulary that captures the relevant phenomena of a specific problem, but still leaves enough generality to be useful to a large class of problems. While this can be said about any mathematical theory, it is especially apparent in the field of nonsmooth differential equations.

\subsection{Nonsmooth mechanical systems with impacts} \label{sect:nonsmoothMechanics}

We have already established that the trajectory $\qe$ of a mechanical system with inequality constraints cannot be the solution of Newton's second axiom~\eqref{eq:newton}, since we must allow $\dot{\qe}$ to have jump discontinuities. To understand the physics behind the system, we have to fall back to more general mechanical principles.

\subsubsection{Hamilton's principle as a differential inclusion}

Classical mechanics offers a few physical principles as governing equations, that are formulated in a variational setting rather than as a second order ordinary differential equation. Hamilton's principle of least action~\eqref{hamstarr} is an example. While there have been several attempts to extend the principle for discontinuous solutions since its development, most of them were undertaken quite recently with the help of the new tools provided by the growing field of nonsmooth dynamical systems, see \cite{erdmann1877,panagiotopoulos2000,fetecau2003,leine2009} and the references therein. It is unclear whether William Hamilton had nonsmooth mechanical systems on his mind in the 1830s. But his principle is formulated with just enough generality to be extended to the nonsmooth case.

To begin, let us consider the problem in a function space setting. While the (generalized) coordinates $\qe$ are not continuously differentiable, they must still be continuous: A material point of a mechanical system cannot just disappear and reappear elsewhere. The velocity $\dot{\qe}$ is continuous almost everywhere, except for a countable number of time points, where the velocity undergoes finite jumps. In mathematical terms, $\qe: \mathbb{R} \supset [0,T] \to \mathbb{R}^N$ is an absolutely continuous function and $\dot{\qe}: [0,T] \to \mathbb{R}^N$ is a function of specially bounded variation,
\[ \qe \in \ac \hskip0.25cm \text{and} \hskip0.25cm \dot{\qe} \in \bv.\]
Introducing inequality constraints $\svec{g}(\qe) \geq 0$ into the system translates to finding a feasible physical solution $\qe$ in the set
\[ \mathcal{S} = \left\{\; \qe \in \ac \; \left| \; \svec{g}(\qe(t)) \geq 0 \hskip0.25cm \forall t \in [0,T]\; \right.\right\}. \]

The classical principle of virtual work and the closely related principle of D'Alem\-bert can be used to examine systems, where the feasible set is a manifold 
\[ \mathcal{M} = \left\{\; \qe \in \ac \; \left| \; \svec{g}(\qe(t)) = 0 \hskip0.25cm \forall t \in [0,T]\; \right.\right\}\]
with a sufficiently smooth function $\svec{g}$. The principle requires the constraint forces $\svec{f}_c(\qe)$ to be orthogonal to the tangent space $T_{\qe}\mathcal{M}$ of the manifold at $\qe$ at all times, i.e. 
\[ -\svec{f}_c(\qe) \in N_{\qe}{\mathcal{M}}, \]
where $N_{\qe}{\mathcal{M}}$ is the annihilator space of $T_{\qe}\mathcal{M}$, see Figure~\ref{fig:tangentspace}. This concept can be extended for the inequality case, $\qe \in \mathcal{S}$, simply by replacing the tangent space with the tangent cone and the annihilator with the normal cone to $\mathcal{S}$ at $\qe$, see~Figure~\ref{fig:tangentcone} \cite{leine2009}:
\begin{equation}
 - \svec{f}_c(\qe) \in N_{\qe}\mathcal{S}.  \label{eq:nonsmoothDalembert}
\end{equation}
Notice, that the normal cone
\[ N_{\qe}\mathcal{S} = \left\{\left.\;  \svec{\xi} \in \ac^* \; \right| \; \langle \svec{\xi}, \svec{d} \rangle \leq 0 \hskip0.25cm \forall \svec{d} \in T_{\qe}\mathcal{S} \; \right\} \]
is a subset of the \emph{dual space} of the space of absolutely continuous functions.

Leine, Aeberhard and Glocker use the nonsmooth version~\eqref{eq:nonsmoothDalembert} of the principle of virtual work and some further mild assumptions\footnote{Particularly, the feasible set $\mathcal{S}$ does not have to be convex, but it is required to be tangentially regular.} to derive a version of Hamilton's principle of least action that is valid in the presence of inequality constraints and discontinuous velocities \cite{leine2009}. The classical version of this principle states that the Fr\'echet--derivative of the action vanishes at $\qe$,
\begin{equation}
 \delta S = \delta \int_0^T L(\qe,\dot{\qe})\D t = 0 \label{eq:hamiltonEq}
\end{equation}
where $L(\qe, \dot{\qe}) = T(\qe,\dot{\qe}) - U(\qe)$ denotes the difference between kinetic and potential energy. In the nonsmooth case, the principle reads
\begin{equation}
 \delta S \in N_{\qe}\mathcal{S}. \label{eq:hamiltonIneq}
\end{equation}
If the feasible set takes up the whole space, $\mathcal{S} = \ac$, the normal cone to $\mathcal{S}$ at any $\qe \in \mathcal{S}$ contains only the zero, $N_{\qe}\mathcal{S} = \{ \svec{0}\}$. In this case, \eqref{eq:hamiltonIneq} reappears as its classical pendant~\eqref{eq:hamiltonEq}.

\subsubsection{Forces and Accelerations are Measures}

It would be convenient to have a form similar to the Lagrange equations of the first kind~\eqref{f1} for inequality constraints that makes the constraint forces explicit. To reach this goal at the end of Section~\ref{sect:lagrangianmult}, we must first accept that the constraint forces of inequality constraints are not necessarily classical functions anymore. If the velocity has a jump, its derivative at that point in time does not exist classically. We can find a weak derivative at the discontinuity point, that is a Dirac delta distribution. This Dirac delta distribution is a measure rather than a function. 

Recall that the constraint forces for inequality constraints are from the dual space of absolutely continuous functions
\[ -\svec{f}_c(\qe) \in N_{\qe} \mathcal{S} \subset (\ac)^*. \]
The space of absolutely continuous functions equipped with the weak norm
\[ \|\qe\| = \max \left\{ \;\sup_{t \in [0,T]}\|\qe(t)\|,\;  \sup_{t \in [0,T]}\|\dot{\qe}(t)\|\;\right\} \]
is a Banach space. Its dual space is given by the signed Radon measures \cite{natanson1977,gerdts2012optimal}, in the sense that any functional $\svec{f} \in (\ac)^*$ can be written as
\[ \langle \svec{f}, \svec{x} \rangle = \int_0^T \svec{x} \;\D \svec{p}, \]
where $\D \svec{p}$ is the Lebesgue--Stieltjes measure of a function of bounded variation~$\svec{p}$.
%
%
It is remarkable how impulsive forces are a direct consequence only of the principle of virtual work applied to the correct function space.

\subsubsection{Existence of Lagrangian multipliers} \label{sect:lagrangianmult}

Up to this point, $\mathcal{S}$ could be any tangentially regular subset of $\ac$. We haven't exploited yet that the set $\mathcal{S}$ is defined by the inequality constraints $\svec{g}(\qe)\geq 0$. In doing so we can characterize the problem using Lagrangian multipliers.

Before the polish mathematician Stanislav Kurcyusz died in a tragic accident at a young age in 1978, he made important contributions to optimization in Banach spaces subject to operator inequality constraints. His theory translates very nicely to the variational formulation of nonsmooth mechanics. Using the results published in \cite{kurcyusz1976,zowe1979} the following theorem can be formulated, that rewrites \eqref{eq:hamiltonIneq} in terms of Lagrangian multipliers \cite{kleinert2015} and a complementarity condition:

\begin{theorem} \label{thm:existlagrange}
 Let $\svec{g}: \mathbb{R}^{n_q} \to \mathbb{R}^{n_\lambda}$ be continuously differentiable and let
 \[\mathcal{S} = \left\{\; \qe \in \ac \; \left| \; \svec{g}(\qe(t)) \geq 0 \hskip0.25cm \forall t \in [0,T]\; \right.\right\} \]
 denote the set of admissible trajectories. Let $\qe \in \mathcal{S}$, let $\mathcal{S}$ be tangentially regular and assume that the \emph{Robinson regularity condition}
 \[ \im (\svec{G}(\qe)) + \spann(\svec{g}(\qe)) - \mathbb{R}_+^{n_\lambda} = \mathbb{R}^{n_\lambda} \]
 holds, where $\im (\svec{G}(\qe))$ denotes the image of the Jacobian of $\svec{g}$ at $\qe$, $\spann (\svec{g}(\qe))$ the space spanned by $\svec{g}(\qe)$ and $\mathbb{R}_+^{n_\lambda}$ the positive orthant in $\mathbb{R}^{n_\lambda}$.
 Then there exists a non--negative measure $\D\svec{\lambda}$, such that
 \begin{subequations}\label{eq:existlagrange}
 \begin{align}
  0 &= \delta S(\delta\qe) + \int_0^T \delta \qe^T \svec{G}(\qe)^T \D \svec{\lambda} \label{eq:existlagrange_1}\\
  0 &= \int_0^T \svec{g}(\qe)^T \D\svec{\lambda}.
 \end{align}
 \end{subequations}
for all variations $\delta\qe \in \ac$ with $\delta \qe(0) = \delta\qe(T)=0$.
\end{theorem}
We can now bring Equation~\eqref{eq:existlagrange_1} in a more recognizable form that resembles the Lagrange equations of first kind~\eqref{f1}. Using $\delta \dot{\qe} \D t = \D (\delta \qe)$, integration by parts for Lebesgue--Stieltjes measures~\eqref{eq:intbyparts} and $\delta \qe(0)= \delta \qe(T)= 0$, we can rewrite the Fr\'echet derivative $\delta S(\delta \qe)$ of the action $S$ as
\begin{align}
 \delta S(\delta \qe) &= \int_0^T \delta \qe^T \frac{\partial L}{\partial \qe} \D t + \int_0^T \delta \dot{\qe}^T \frac{\partial L}{\partial \dot{\qe}} \D t \notag \\
 &= \int_0^T \delta \qe^T \frac{\partial L}{\partial \qe} \D t + \int_0^T \frac{\partial L}{\partial \dot{\qe}}^T \D (\delta \qe) \notag \\
 &= \int_0^T \delta \qe^T \frac{\partial L}{\partial \qe} \D t + \underbrace{\left[ \delta \qe^T \frac{\partial L}{\partial \dot{\qe}} \right]_0^T}_{=0} - \int_0^T \delta \qe^T \D \left( \frac{\partial L}{\partial \dot{\qe}} \right) \notag \\
 &=\int_0^T \delta \qe^T \frac{\partial L}{\partial \qe} \D t - \int_0^T \delta \qe^T \D \left( \frac{\partial L}{\partial \dot{\qe}} \right). \label{eq:deltaS}
\end{align}
Plugging~\eqref{eq:deltaS} into~\eqref{eq:existlagrange_1} yields
\begin{equation}
 0= \int_0^T \delta \qe^T \frac{\partial L}{\partial \qe} \D t - \int_0^T \delta \qe^T \D \left( \frac{\partial L}{\partial \dot{\qe}} \right) + \int_0^T \delta \qe^T \svec{G}(\qe)^T\D \svec{\lambda} \label{eq:variational-eulerlagrange}
\end{equation}
or, equivalently expressed as an equality of measures,
\begin{equation}
 \D \left( \frac{\partial L}{\partial \dot{\qe}} \right) = \frac{\partial L}{\partial \qe} \D t + \svec{G}(\qe)^T\D \svec{\lambda}. \label{eq:mde-eulerlagrange}
\end{equation}
These are the nonsmooth version of the Euler--Lagrange equations. As long as the momentum 
$\partial L/\partial \dot{\qe}$
has discontinuous jumps, i.e. is not absolutely continuous, it holds
\[ \D \left( \frac{\partial L}{\partial \dot{\qe}} \right) \neq \frac{\D}{\D t} \frac{\partial L}{\partial \dot{\qe}} \; \D t.
\]
Therefore we cannot derive an ordinary differential equation from the variational problem~\eqref{eq:variational-eulerlagrange} using the fundamental lemma of calculus, as is usually done to arrive at the Euler--Lagrange equations. The physics can only be described variationally, in this case in the form of a measure differential equation~\eqref{eq:mde-eulerlagrange}.

Finally, by plugging in the Lagrange function of a multibody system
\[ L(\qe,\dot{\qe}) = \frac{1}{2}\dot{\qe}^T\svec{M}(\qe)\dot{\qe} - U(\qe) \]
and using
\begin{align}
 \frac{\partial L }{\partial \qe} &= \frac{1}{2}\dot{\qe}^T\frac{\partial \svec{M}(\qe)}{\partial \qe}\dot{\qe} - \frac{\partial U(\qe)}{\partial \qe}, \notag \\
 \D\left( \frac{\partial L }{\partial \dot{\qe}} \right) &= \svec{M}(\qe)\D \dot{\qe} + \dot{\qe}^T \D \svec{M}(\qe) \notag \\
 &=\svec{M}(\qe)\D \dot{\qe} + \dot{\qe}^T \dot{\svec{M}}(\qe)\D t \notag \\
 &=\svec{M}(\qe)\D \dot{\qe} + \dot{\qe}^T \frac{\partial \svec{M}(\qe)}{\partial \qe} \dot{\qe} \D t,\notag \\
 \svec{f}(\qe,\dot{\qe},t) :&= -\frac{1}{2}\dot{\qe}^T\frac{\partial \svec{M}(\qe)}{\partial \qe}\dot{\qe} - \frac{\partial U(\qe)}{\partial \qe}, \notag
\end{align}
Equation~\eqref{eq:mde-eulerlagrange} becomes
\begin{equation}
 \svec{M}(\qe)\;\D \dot{\qe} = \svec{f}(\qe,\dot{\qe},t)\;\D t + \svec{G}(\qe)^T\;\D \svec{\lambda}. \label{eq:mde-mbd}
\end{equation}
Replacing Equation~\eqref{eq:existlagrange_1} in Theorem~\ref{thm:existlagrange} with Equation~\eqref{eq:mde-mbd} yields the Lagrange equations of first kind for the nonsmooth case with impacts
\begin{subequations}\label{eq:f2}
\begin{align}
 \svec{M}(\qe)\;\D \dot{\qe} &= \svec{f}(\qe,\dot{\qe},t)\;\D t + \svec{G}(\qe)^T\;\D \svec{\lambda}, \label{eq:f2-a}  \\
 \svec{0} &\leq \svec{g}(\qe), \label{eq:f2-b} \\
 \svec{0} &\leq \D \svec{\lambda}, \label{eq:f2-c} \\
 0 &= \int_0^T \svec{g}(\qe)^T \;\D \svec{\lambda}. \label{eq:f2-d}
\end{align}
\end{subequations}
The differential equation~\eqref{f1a} in the smooth case is replaced by the MDE~\eqref{eq:f2-a} in the nonsmooth version and the equality constraint~\eqref{f1b} turns into the complementarity conditions~\eqref{eq:f2-b}--\eqref{eq:f2-d}. Due to the positivity \eqref{eq:f2-b} and \eqref{eq:f2-c}, \eqref{eq:f2-d} means that the measure $\D \svec{\lambda}$ is nonzero only for subsets of $[0,T]$ on which $\svec{g}$ is zero.

\subsection{Numerical solution strategies} \label{sect:nonsmoothNumerics}




Nonsmooth dynamical systems have been researched rigorously since the late 1950s and many theoretical results were developed during the 1970s and 1980s. However, based on the number of publications on the subject, the development of numerical methods only slowly took up pace in the 1990s and led to a boom in the field just recently during the 2000s and 2010s.

In the 1990s, DAEs were a hot topic and many multibody simulation tools were developed. Contacts in multibody dynamical systems were usually modeled by penalizing interpenetration of rigid bodies with a stiff spring--damper element at the contact point. While this is a simple and useful model for many applications, it can be somewhat restrictive for mechanical systems with many contacts. At least two additional parameters per contact, a stiffness and a damping coefficient, have an influence on the overall behavior of the model and must be chosen appropriately. Small time steps must be used, so that a switch in contact state is not missed. Very stiff spring--damper elements can only be used together with very small time steps. As more computational power became available, industrial scale simulations of granular material became viable. The classical Discrete Element Method of Cundall and Strack \cite{Cundall1979} is based on a penalized contact model. Very small time steps must be used to maintain a stable simulation. Nonsmooth dynamical systems, where contacts are modeled using inequality constraints, became an attractive alternative to deal with contacts and collisions in multibody simulations. With the new approach, it is not necessary to capture the change in contact states exactly and the ``infinitely stiff'' character of the contact laws do not yield stiff differential systems. Much larger time steps can be used.

\subsubsection{Even--driven and even--capturing methods}

There is a very intuitive method for dealing with nonsmooth events when solving the equations of motion of a dynamical system. Remember, that the velocity is a function of bounded variation, and as such is smooth almost everywhere except at a countable number of discontinuity points. \emph{Event--driven} integrators calculate the smooth trajectory of the system using available ODE or DAE solvers, until such an event is detected. The time integration stops, the event can be handled, e.g. by evaluating a contact model, and the integrator is restarted with new initial conditions until the next nonsmooth event is detected. Today, many ODE and DAE solvers integrate event--detection features. These techniques are useful if the overall number of nonsmooth events is not too large and the time points of the events can be predicted effectively.

But already the simple bouncing ball problem with a restitution coefficient $0 < e < 1$ displays Zeno behavior, that is it has an accumulation point of discontinuous events in time. A rigorous event--driven time integrator would not be able to pass this accumulation point without any additional trickery. In addition to this, there are numerous applications where the frequency of nonsmooth events in time is high and the numerical overhead of event--detection becomes the bottleneck of the entire simulation. This is especially true for dynamic simulations of granular matter, where each rigid particle is subject to frequent changes in the contact state. For this reason, this section focuses solely on time integration schemes that can step over nonsmooth events and capture the net movement during time steps, regardless of whether this movement is due to smooth motion or a nonsmooth phenomenon. Integration schemes of this type are often denoted as \emph{event--capturing} or \emph{time--stepping} methods.

Time--Stepping methods are commonly split into two categories, namely the Paoli--Schatzman \cite{paoli2002,paoli2002a} and the Moreau--Jean time--stepping schemes \cite{jourdain1998,Jean1999}. The main difference between the two is that the first aims at finding a solution on position level while the latter tackle the problem on a velocity level. In both methods projections onto certain sets must be performed. In the first, projections onto arbitrary convex sets are applied, while in the latter only projection onto the normal and tangent cones to these sets play a role, which are easier to compute in general \cite{acary2008}. In the following, we will concentrate the discussion on Moreau--Jean time--stepping.

\subsubsection{Nonsmooth time--stepping}

Consider a discrete version of Equation~\eqref{eq:f2-a}, that can be obtained by integrating over the time interval $[t_k, t_{k+1}]$:
\begin{equation}
 \svec{M}(\qe_k) \left( \svec{v}_{k+1} - \svec{v}_k \right) =  \svec{k}_{k} + \svec{G}(\qe_k) \svec{p}_{k+1} \label{eq:timestepping-v}
\end{equation}
where $h =t_{k+1}-t_k$, $\qe_k \approx \qe(t_k)$ and $\svec{v}_k \approx \svec{v}(t_k):= \dot{\qe}(t_k)$,
\[\svec{p}_{k+1} = \int_{t_k}^{t_{k+1}}\D \svec{\lambda}\]
is an impulse that appears as a new unknown and
\[ \svec{k}_k \approx \int_{t_k}^{t_{k+1}}\svec{f}(\qe(\tau),{\svec{v}}(\tau),\tau)\;\D\tau\]
is an approximation of the force integral that can be obtained using a suitable quadrature rule.  The generalized position of the mechanical system can be calculated for instance using a $\Theta$--method
\begin{equation}
 \qe_{k+1} = \qe_k + h ( \Theta \svec{v}_{k+1} + (1-\Theta)\svec{v}_k ). \label{eq:timestepping-q}
\end{equation}
With $\Theta = 0$, Equation~\eqref{eq:timestepping-q} corresponds to the explicit Euler, $\Theta=\frac{1}{2}$ yields the mid--point rule and $\Theta = 1$ the implicit Euler method.

For simplicity, it is assumed here, that the mass matrix $\svec{M}$ and constraint Jacobian $\svec{G}$ remain almost constant during a time step and that $\svec{f}(\qe(\tau),{\svec{v}}(\tau))$ can be approximated or extrapolated for the time interval $[t_k, t_{k+1}]$. For applications in linear elasticity, $\svec{f}$ is often linear in $\qe$ and $\svec{v}$, e.g.
\[ \svec{f} = \svec{f}_{\text{ext}} - \svec{K} \qe - \svec{D} \svec{v}, \]
with stiffness and damping matrices $\svec{K}$ and $\svec{D}$. If in this case the $\Theta$--method is used for the approximation of the force integral $\svec{k}_k$, no further extrapolation is needed, see for instance \cite{acary2008}. Convergence results for this time--stepping scheme are provided for instance in \cite{monteiromarques1993,stewart1998,acary2008}. 

In \cite{schindler2012,Schindler2014,kleinert2017} a more rigorous derivation of the discrete time--stepping equations is provided. Similarly to the way finite element methods are constructed from a weak form of a partial differential equation, the authors consider equation~\eqref{eq:f2} as a variational problem and apply a discontinuous Galerkin approximation in time to obtain the above mentioned time stepping equations. 

The Lagrangian multiplier appears as as an impulse $\svec{p}_{k+1}$, that has the unit of momentum. No approximations for the accelerations and forces exist. The dynamics are solved directly on velocity level. The new unknown $\svec{p}_{k+1}$ does not distinguish between smooth or nonsmooth parts of the acceleration. This can potentially be a problem in applications, where particularly the loads due to contact forces are of interest \cite{kleinert2013}.

In order to advance the time--stepping scheme, we need to solve for $\svec{p}_{k+1}$. The complementarity problem~\eqref{eq:f2-b} -- \eqref{eq:f2-c} for $\D \svec{\lambda}$ translates directly to a complementarity problem for the new unknown $\svec{p}_{k+1}$.
\begin{equation}
 \svec{g}(\qe_{k+1}) \geq \svec{0}, \hskip0.5cm \svec{p}_{k+1} \geq 0, \hskip0.5cm \svec{p}_{k+1}^T\svec{g}(\qe_{k+1}) = 0. \label{eq:timestepping-ncp}
\end{equation}
In general, this is a \emph{nonlinear complementarity problem} (NCP)\index{complementarity problem!nonlinear}. For $\Theta=1$, Equations~\eqref{eq:timestepping-v},\eqref{eq:timestepping-q} and the NCP~\eqref{eq:timestepping-ncp} amount to a nonsmooth version of the SHAKE integrator~\eqref{shake} \cite{HaLW02,hairer2003}. The position $\qe_{k+1}$ is updated using the velocity $\svec{v}_{k+1}$, which is calculated in such a way, that the constraint is satisfied at the end of the time step from $t_k$ to $t_{k+1}$.

\begin{figure}
 \centering
 
 \def\sim{100}
 \def\xlimlo{0}
 \def\xlimhi{30}
 \def\vlimlo{-20}
 \def\vlimhi{2}
 
 \subcaptionbox{$t=0.0\;\mathrm{s}$\label{fig:snap1}}
 {\centering\includegraphics[width=0.47\columnwidth]{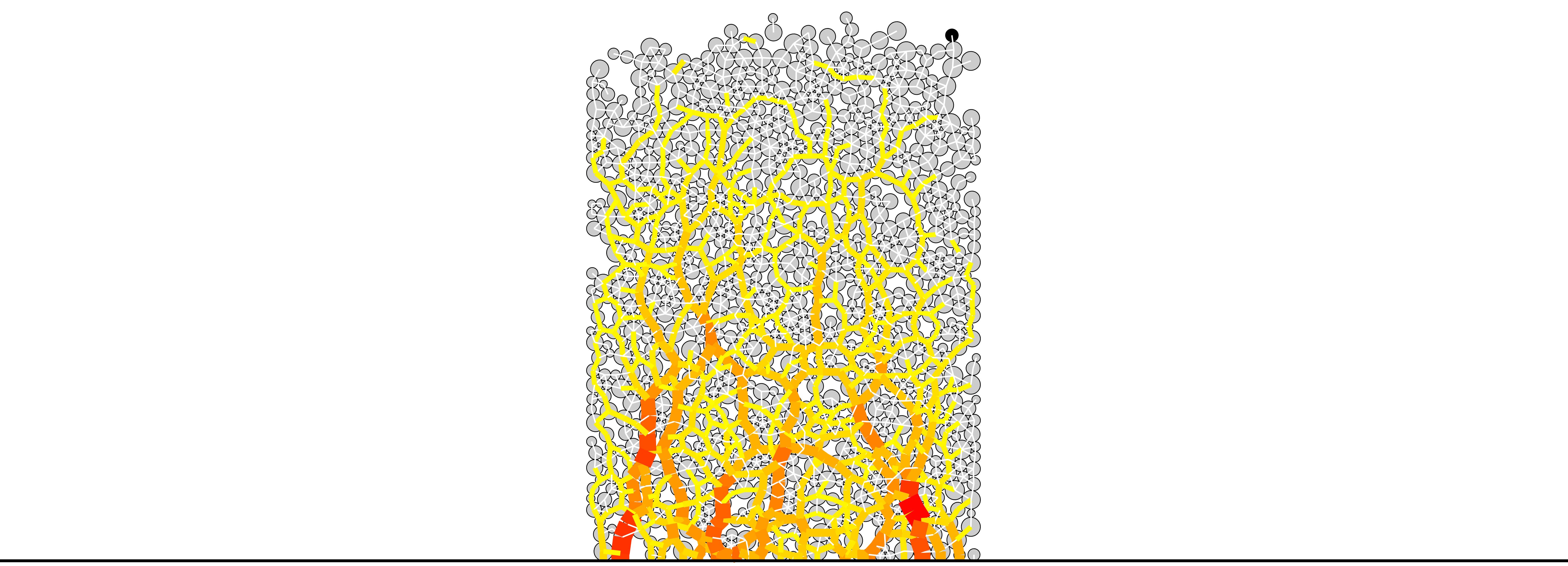}}\hfill
 \subcaptionbox{$t=2.0\;\mathrm{s}$\label{fig:snap2}}
 {\centering\includegraphics[width=0.47\columnwidth]{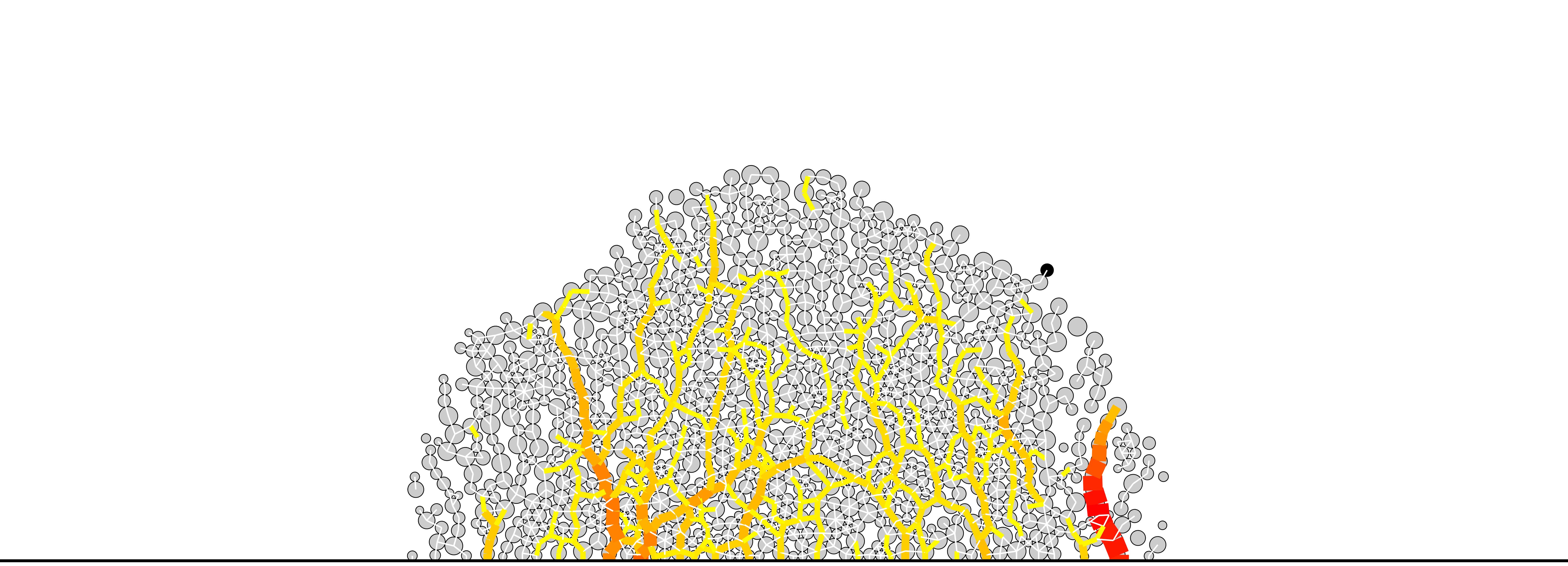}}
 
 \subcaptionbox{$t=3.0\;\mathrm{s}$\label{fig:snap3}}
 {\centering\includegraphics[width=0.47\columnwidth]{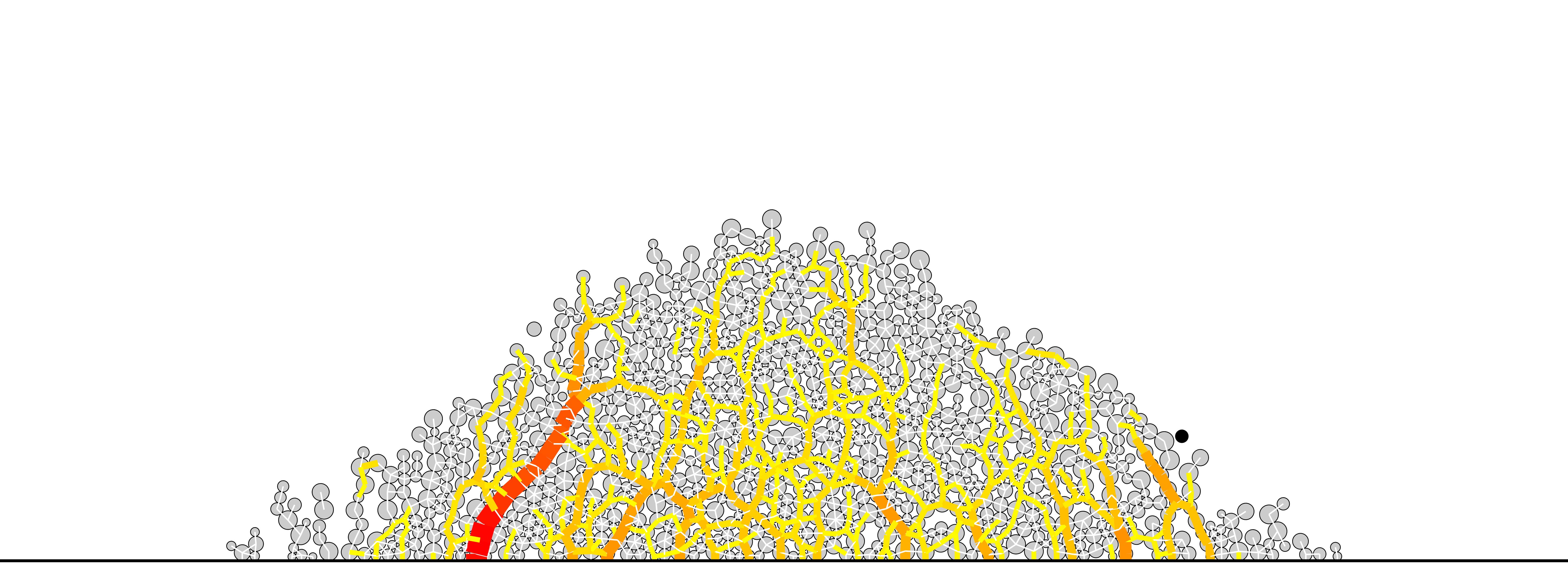}}\hfill
 \subcaptionbox{$t=5.0\;\mathrm{s}$\label{fig:snap4}}
 {\centering\includegraphics[width=0.47\columnwidth]{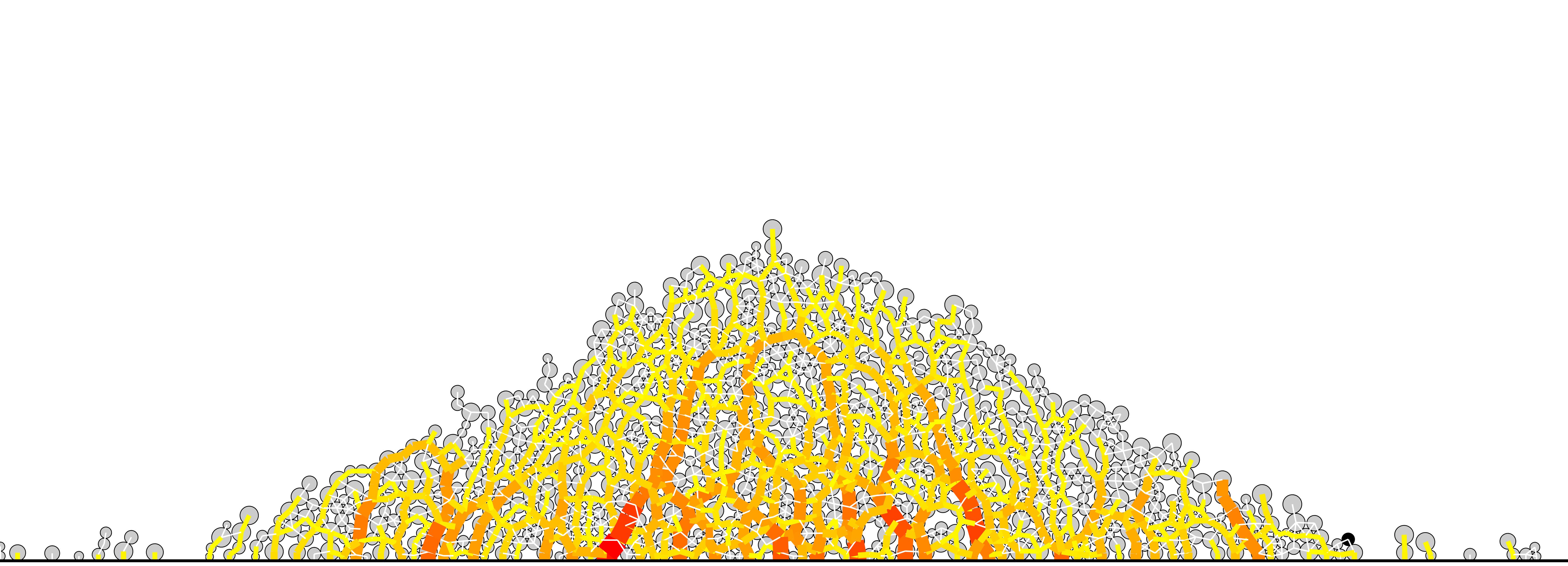}}
 
 \subcaptionbox{The continuous, piecewise linear position and discontinuous, piecewise constant velocity of a single particle in vertical direction.\label{fig:pilexv}}{
 \centering\begin{tikzpicture}
 
 \pgfplotstableread[col sep = comma]{data_TXV_\sim.dat}\mydata;
 
  \begin{axis}[
  width=0.8\columnwidth,
  height=0.4\columnwidth,
  axis y line*=left,
  xmin=0, xmax=5,
  ymin=\xlimlo, ymax=\xlimhi,
  xlabel={$t$ [s]},
  ylabel={$x_2$ [mm]}]
  
\addplot [ 
           color=blue,
           style=dashed
         ]
         table[ x index=0,y index=2] {\mydata};
\label{plot_one}
\end{axis}

\begin{axis}[
  width=0.8\columnwidth,
  height=0.4*\columnwidth,
  axis y line*=right,
  xmin=0, xmax=5,
  ymin=\vlimlo, ymax=\vlimhi,
  axis x line=none,
  ylabel={$v_2$ [mm/s]},
  legend style = {at={(0.05,0.05)},anchor=south west}]
\addlegendimage{/pgfplots/refstyle=plot_one}\addlegendentry{position $x_2$}
\addplot[jump mark left,red, line width=0.5pt] 
  table[ x index=0,y index=4] {\mydata};
\addlegendentry{velocity $v_2$}
\end{axis}

 \end{tikzpicture}
 }
 
 \caption{A simulation of a collapsing pile of 1000 rotationless disks in two dimensions. Figures~\ref{fig:snap1}--\ref{fig:snap4} show snapshots at different time points in the simulation. The images contain visualizations of the so--called \emph{force chains} in the granular material. The vertical position and velocity of the highlighted particle is shown in Figure~\ref{fig:pilexv}.}
 \label{fig:pile-collapse}
 
\end{figure}
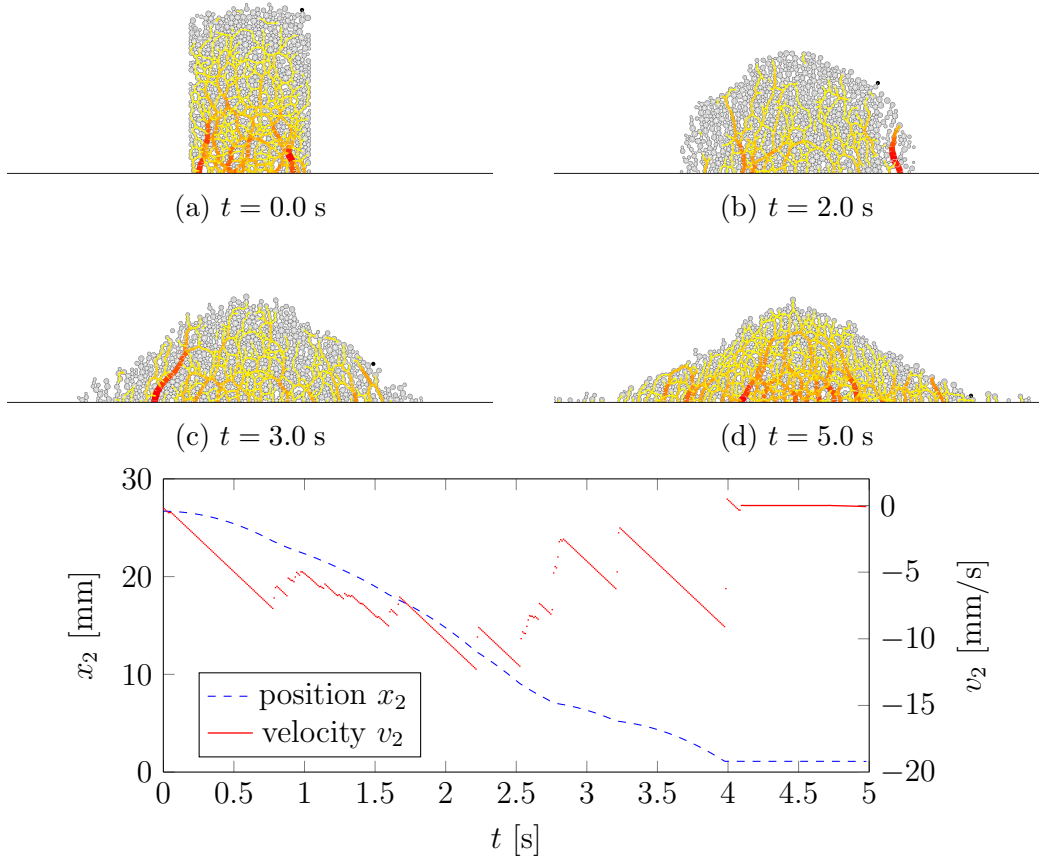

Figure~\ref{fig:pile-collapse} shows a simulation of two--dimensional rotationless circular particles. In this simulation, the SHAKE integrator was used. The velocity is piecewise constant and the position is continuous and piecewise linear.

\subsubsection{Dealing with collisions}

How can we make sure, that the constraints are satisfied at the end of the time step, if we do not use any event detection in our numerical method? There are two common strategies to do this at least approximately. 

In the first strategy, we allow the constraints to be violated and assume that the time step size is small enough, so that these violations remain small. Let $\mathcal{I} = \left\{ i_1, ..., i_n \in \mathbb{N} \right\}$ denote the index set of the active and the violated inequality constraints at time $t_k$,
\[ g_{s}(\qe_k) \leq 0, \hskip0.5cm s \in \mathcal{I}. \]
For this \emph{active set} of constraints it it clear, that the contact velocity must be non--negative by the next time step, so that the constraint is not violated any further,
\begin{equation}
 \dot{g_{s}}(\qe_{k+1}) = \frac{\partial g_{s}(\qe_{k+1})}{\partial \qe}^T \cdot  \svec{v}_{k+1} \geq 0 \hskip0.5cm s \in \mathcal{I}. \label{eq:constraint-v1}
\end{equation}
Thus, with the additional assumption that the constraint Jacobian
\[ \svec{G}_{k,s} = \frac{\partial g_{s}(\qe_{k})}{\partial \qe} \approx \frac{\partial g_{s}(\qe_{k+1})}{\partial \qe} \]
remains almost constant during one time step we can replace the NCP~\eqref{eq:timestepping-ncp} with the following complementarity problem on velocity level, that is solved only for the active set of inequality constraints
\[ \svec{G}_{k,s} \cdot  \svec{v}_{k+1} \geq 0, \hskip0.5cm p_{k+1,s} \geq 0, \hskip0.5cm p_{k+1,s}\svec{G}_{k,s} \cdot  \svec{v}_{k+1} = 0 \hskip0.25cm\forall\;s\in\mathcal{I}. \]
Here $p_{k+1,s}$ denotes the component of $\svec{p}_{k+1}$ associated to the inequality constraint with index $s$. Note that initial violations remain, as the reaction impulse $p_{k+1,s}$ is just large enough to not decrease the value of $g_s$ any further.

The other strategy consists of linearizing the inequality constraint 
around $\qe_k$,
\begin{align}
 \svec{g}(\qe_{k+1}) &= \svec{g}\left(\qe_k + h \left( \Theta \svec{v}_{k+1} + (1-\Theta)\svec{v}_k \right) \right)   \notag \\
 &\approx \svec{g}(\qe_k) + h \frac{\partial \svec{g}(\qe_{k})}{\partial \qe} ^T\left( \Theta \svec{v}_{k+1} + (1-\Theta)\svec{v}_k \right) \geq 0. \label{eq:timestepping-linearized}
\end{align}
In this strategy we are looking ahead: We can calculate the velocity $\svec{v}_{k+1}$ in such a way, that the inequality constraint is satisfied at the end of the next time step. We do not need to separate the constraints into an active and inactive set. If the value of a linearized inequality constraint is positive and far away from zero, the NCP~\eqref{eq:timestepping-ncp} implies that he corresponding contact impulse must be zero. Dividing the linearized inequality constraint~\eqref{eq:timestepping-linearized} by the time step size $h$ reveals it as an inequality constraint on velocity level with an additional stabilization term \cite{anitescu2004b}
\begin{equation}
 \svec{u}_{k+1} = \frac{\svec{g}(\qe)}{h} + \frac{\partial \svec{g}(\qe_{k})}{\partial \qe}^T \left( \Theta \svec{v}_{k+1} + (1-\Theta)\svec{v}_k \right) \geq 0, \label{eq:constraint-v2}
\end{equation}
very similar to Baumgarte stabilization of velocity based bilateral constraints \cite{Baum72}. 

The NCP~\eqref{eq:timestepping-ncp} is replaced by
\[ \svec{u}_{k+1} \geq \svec{0}, \hskip0.5cm \svec{p}_{k+1} \geq 0, \hskip0.5cm \svec{p}_{k+1}^T\svec{u}_{k+1} = 0. \]
As the constraints are considered on velocity level, these two strategies are very similar to index reduction methods for DAEs.

Both approaches amount to perfectly inelastic collisions. The unilateral constraint can be modified so that partially or perfectly elastic collisions are modeled, e.g. by incorporating Newton's impact law. In general however, we must separate the collision into a compression and decompression phase to accomplish this without introducing additional errors. This requires the solution of two complementarity problems per time step \cite{Anitescu1997a}. And even then, this only models elastic collisions between two bodies. modeling inelastic shocks through a network of contacting rigid bodies, such as for example Newton's cradle, is still a challenge in event--capturing strategies \cite{nguyen2018}.

\subsubsection{Dealing with complementarity}

Since the early days of numerical methods for nonsmooth dynamical systems and regardless of which concept from Section~\ref{sect:dviHistory} is used for the description of the problem, most authors formulate the discrete numerical method in terms of a complementarity problem similar to the one from the previous section \cite{moreau1966, Loe82, loetstedt1984}. Posing the inequality constraints~\eqref{eq:constraint-v1} or~\eqref{eq:constraint-v2} on velocity level has a big advantage: they are linear in the velocity $\svec{v}_{k+1}$ and therefore also in the unknown impulse $\svec{p}_{k+1}$. Solving~\eqref{eq:timestepping-v} for $\svec{v}_{k+1}$ and inserting the result into~\eqref{eq:constraint-v2} reveals this,
\begin{align}
 \svec{u}_{k+1} =& \frac{\svec{g}(\qe_k)}{h} + \svec{G}_k^T \left( \Theta \svec{v}_{k+1} + (1-\Theta)\svec{v}_k \right) \notag \\
 =& \frac{\svec{g}(\qe_k)}{h} + \svec{G}_k^T \left( \svec{v}_k + \Theta \svec{M}_k^{-1} \svec{k}_k \right) \notag \\
 &+ \Theta \svec{G}_k^T\svec{M}_k^{-1}\svec{G}_k \svec{p}_{k+1} \notag \\
 =& \svec{b}_k + \svec{A}_k \svec{p}_{k+1}, \notag 
\end{align}
where $\svec{G}_k = \svec{G}(\qe_k)$ is the constraint Jacobian, $\svec{M}_k = \svec{M}(\qe_k)$ is the mass matrix, 
\[ \svec{b}_k = \frac{\svec{g}(\qe_k)}{h} + \svec{G}_k^T \left( \svec{v}_k + \Theta \svec{M}_k^{-1} \svec{k}_k \right) \]
and 
\[ \svec{A}_k = \Theta \svec{G}_k^T\svec{M}^{-1}\svec{G}_k \]
is the \emph{Delassus} matrix. The NCP~\eqref{eq:timestepping-ncp} is transformed to a \emph{linear complementarity problem} (LCP)\index{complementarity problem!linear}
\begin{equation}
 \svec{u}_{k+1} = \svec{A}_k \svec{p}_{k+1} + \svec{b}_k \geq 0, \hskip0.5cm \svec{p}_{k+1} \geq 0, \hskip0.5cm \svec{p}_{k+1}^T\svec{u}_{k+1} = 0,\label{eq:timestepping-lcp}
\end{equation}
which is equivalent to the \emph{quadratic program} (QP)\index{quadratic program}
\begin{align}
\begin{split}
 &\min_{\svec{p}_{k+1}}\; \svec{p}_{k+1}^T \svec{A}_k \svec{p}_{k+1}  + \svec{p}_{k+1}^T\svec{b}_k \\
 & \begin{array}{ll} 
 \text{s.t.} & 0 \leq \svec{A}_k \svec{p}_{k+1} + \svec{b}_k,  \\
 & 0 \leq \svec{p}_{k+1}.
   \end{array}
\end{split}
\label{eq:qp}
\end{align}
These LCPs or QPs were usually solved using Lemke's algorithm \cite{murty88}, see for example \cite{loetstedt1984, baraff1992, stewart1996,Anitescu1997a,anitescu2004}. Lemke's algorithm a pivoting strategy for LCPs similar in spirit to Gau{\ss}ian elimination. Lemke's algorithm also has a similar complexity bound as Gau{\ss}ian elimination, which makes it impractical as soon as the complementarity problem becomes quite large. This is the case when many bodies are connected through a network of contacts, as for instance in granular assemblies.

\subsubsection{Coulomb friction and the Painlev\'e paradox}
%

Coulomb friction can be formulated as a nonlinear complementarity problem, see for instance \cite{1998DeSaxce}. The nonlinearity can be removed by approximating the Coulomb friction cone by a polyhedral shape. The NCP for friction turns into an LCP that can be solved together with the LCP for unilateral contacts using Lemke's algorithm. Within this polyhedral approximation, the frictional force can only be exerted along a finite number of predefined directions that fan out from the contact point. This has two major disadvantages. Firstly, we introduce artificial anisotropy if we choose coarse approximations of the Coulomb friction cone. Secondly, for every possible direction of the frictional force a new unknown is introduced. As a consequence, for fine approximations of the Coulomb friction cone, the LCP becomes very large. 

David Stewart uses a convergence proof of the numerical methods published in \cite{stewart1996,Anitescu1997a} that use a polyhedral approximation of the friction cone to resolve the much discussed Painlev\'e paradox \cite{stewart1998}. A detailed analysis of the paradox, first introduced in 1895, is provided in \cite{champneys2016}. It consists of a simple example of a rigid body in frictional contact with a flat rigid surface. As the rigid body slides over the surface, there are instances in time where the system seems to have infinite solutions and instances where it does not seem to have any solutions. This seeming inconsistency of rigid body dynamics with Coulomb friction can be resolved by allowing shocks, i.e. impulsive forces, even when there are no impacts in the system. Stewart shows that the numerical methods of \cite{stewart1996,Anitescu1997a} applied to the Painlev\'e example converge to an impulsive solution of a measure differential inclusion that describes rigid body dynamics subject to Coulomb friction. 

\subsubsection{Augmented Lagrangian and projected Gau{\ss}--Seidel}\index{augmented Lagrangian method}

It has already been established that the complexity bound for Lemke's algorithm is too restrictive for large problems. The publication \cite{landers1986} includes an alternative approach to deal with frictionless contact problems in the context of finite element methods based on an augmented Lagrangian method. This ideas was later adapted by Alart and Curnier in \cite{Alart91} to frictional contact problems. The authors formulate a variational minimization problem at each time step. A energy functional $\varphi(\qe)$ must be minimized by a trajectory $\qe$ without violating inequality constraints $g_i(\qe) \geq 0$, $i=1,...,n_{\lambda}$. This variational inequality--constrained minimization problem is solved using the augmented Lagrangian method. The basic idea is to recast this as an unconstrained optimization problem. For the frictionless case, this leads to
\begin{equation}
 \min \Phi_r(\qe,\svec{\gamma}) =  \varphi(\qe) - \frac{1}{2r}\sum_{i=1}^{n_\lambda}  \left( \| \gamma_i \| ^2 - \dist^2(\gamma_i + r g_i(\qe),\mathbb{R}_+) \right),  \label{eq:al}
\end{equation}
where
\[ \dist(x,C) = x - \proj_C(x) \]
denotes the Euclidean distance of a point $x$ from a set $C$ and $r>0$ is a numerical factor controlling the steepness of $\Phi_r$. The complementarity conditions and inequality constraints are eliminated from~\eqref{eq:al}. The unknown $\svec{\gamma}$ is an approximation of the Lagrangian multiplier associated with the inequality--constrained variational problem. A solution to the minimization problem~\eqref{eq:al} is sought by solving the saddle point problem $\nabla \Phi_r(\qe,\svec{\gamma})=0$. In the context of time--stepping methods, the saddle point problem is solved together with~\eqref{eq:timestepping-v}, where $\svec{\gamma}$ takes on the role of the reaction impulse $\svec{p}_{k+1}$. Alart and Curnier suggest a generalized Newton method for the solution, that boils down to the following iterative projection method \cite{Alart91,leine2004,Pfeiffer2006,studer2008}. 
Let $i=0$. Until convergence is achieved in $\svec{p}_{k+1}^{(i)}$, repeat the following steps:
\begin{enumerate}
 \item $\svec{v}_{k+1} = \svec{v}_k + \svec{M}_k^{-1}(\svec{k}_k + \svec{G}_k \svec{p}_{k+1}^{(i)}),$ where $\svec{p}_{k+1}^{(i)}$ is given from a previous iteration or from an initial guess.
 \item Calculate $\svec{u}_{k+1}$ from $\svec{v}_{k+1}$ using \eqref{eq:constraint-v2}.
 \item $\svec{p}_{k+1}^{(i+1)} = \proj_{\mathbb{R}^{n_{\lambda}}_+}(\svec{p}_{k+1}^{(i)} - r\svec{u}_{k+1})$.
 \item $i \gets i+1$.
\end{enumerate}
Note, that for any $r > 0$ the following equivalence holds
\begin{equation}
 \svec{p} = \proj_{\mathbb{R}_+^{n_{\lambda}}}(\svec{p}-r\svec{u}) \hskip0.5cm \Leftrightarrow \hskip0.5cm \svec{p} \geq \svec{0}, \hskip0.25cm \svec{u} \geq \svec{0}, \hskip0.25cm \svec{p}^T\svec{u} = 0. \label{eq:pgj}
\end{equation}
In other words, if the iteration converges, it converges to a solution of the complementarity problem~\eqref{eq:timestepping-lcp}.
For brevity, friction has not been considered here, details can be found in the before mentioned publications. Equation~\eqref{eq:pgj} is written as a vector equation.
We can rewrite the projection step as
\[ \svec{p}_{k+1}^{(i+1)} = \proj_{\mathbb{R}^{n_{\lambda}}_+}\left(\svec{p}_{k+1}^{(i)} - r\left( \svec{A}_k \svec{p}_{k+1}^{(i)} + \svec{b}_k \right)\right), \]
which reveals it as the Gau{\ss}--Jacobi method for the solution of linear systems of equations with an additional projection, called the \emph{projected Gau\ss--Jacobi method method} (PGJ)\index{projected Gau"s--Jacobi method}. Of course, given the previous iterate $\svec{p}_{k+1}^{(i)}$, the equation can be evaluated component--wise and the newly calculated components $j=1,...,l < n_{\lambda}$ can directly be used in the right--hand side of the equation for the calculation for the component with index $l+1$. This method is called the \emph{projected Gau\ss--Seidel method} (PGS)\index{projected Gau\ss--Seidel method}. It has slightly better convergence properties than PGJ, but PGJ lends itself for a parallel implementation. The convergence of both PGJ and PGS can be controlled with a \emph{successive overrelaxation} parameter $\alpha$,
\begin{align}
 \tilde{\svec{p}}_{k+1}^{(i+1)} &= \proj_{\mathbb{R}^{n_{\lambda}}_+}\left(\svec{p}_{k+1}^{(i)} - r\left( \svec{A}_k \svec{p}_{k+1}^{(i)} + \svec{b}_k \right)\right) \notag \\ & \\
 \svec{p}_{k+1}^{(i+1)} &= \alpha \cdot \tilde{\svec{p}}_{k+1}^{(i+1)} + (1-\alpha)\cdot \svec{p}_{k+1}^{(i)}, \notag
\end{align}
yielding the class of \emph{projected successive overrelaxation} (PSOR)\index{projected successive overrelaxation} methods. 

The augmented Lagrangian method is applied to granular simulations in \cite{fortin2002}. In \cite{tasora2008,anitescu2010,tasora2011} a numerical method based on the PGJ, PGS and PSOR schemes is proposed for general cone complementarity problems. The authors use this method to solve for the normal and tangential reaction impulses involved in Coulomb friction simultaneously by directly projecting the three-dimensional contact force onto the Coulomb friction cone. They demonstrate a large number of numerical examples including the simulation of large granular assemblies and provide implementation details in a HPC context.

\subsubsection{Recent developments}

Over the course of the last two decades several simulation codes were developed around nonsmooth mechanical systems. It is not our intention to compile a complete list, but a few of these codes deserve to be mentioned.

Jean and Dubois initiated the open source software LMGC90 \cite{jourdain1998,dubois2011,lmgc_gitlab} which is currently being developed at the University of Montpellier. SICONOS \cite{siconos_github} is an open source software developed at INRIA in Grenoble by the TRITOP team led by Vincent Acary, following the previous works of Bernard Brogliato's BIPOP team. DynamY \cite{dynamy_online} is another example for a C++ library for nonsmooth mechanical systems that was developed during the course of the PhD thesis \cite{studer2008} at the ETH Z\"urich. Algoryx \cite{algoryx_online} is a spin--off software company from Ume{\aa} University, that develops commercial software for the simulation of nonsmooth mechanical systems. Another free simulation code is the PE Rigid Body Physics Engine \cite{preclik2014,pe_online} developed at the University of Erlangen. Finally, the open source software Chrono \cite{chrono_github} is currently being developed by a large team at the University of Wisconsin and the Universit\`a di Parma by Tasora Negrut, Serban and a team of their students.

Recent research in the field of nonsmooth mechanical systems concentrates either on solving large complementarity problems in a short time, or on increasing the order of the integration methods. The first is useful for applications with many rigid bodies in close contact, as is the case for the simulation of granular material. The latter is interesting for general flexible multibody simulation scenarios with a small number of contacts. Here, the presence of unilateral contact and friction decreases the overall order of the integration methods compared to standard DAE solvers, if the Moreau--Jean or Paoli--Schatzman integrators must be used to resolve the collisions consistently.

\paragraph{Fast solvers for complementarity problems}

So far, only two strategies have been discussed to deal with the complementarity conditions. The first was the direct solution of an LCP using a pivoting strategy such as Lemke's method. The other was an augmented Lagrangian approach that uses a fix point iteration for the Lagrangian multipliers, which involves a projection onto a feasible set.

Other numerical tools from the field of continuous optimization than Lemke's method can be borrowed to solve the complementarity problems in multibody simulations. For this purpose, it makes sense to introduce two categories of iterative solvers. Consider for this the QP~\eqref{eq:qp}, where the objective function takes on the form
\[ f(\svec{p}) = \svec{p}^T\svec{A} \svec{p} + \svec{p}^T\svec{b}.\]
The first category consists of those methods, that only require  the gradient $\nabla f(\svec{p}) = \svec{A} \svec{p} + \svec{b}$
in every iteration. The second category consists of those methods, that need to solve linear systems involving the Hessian $\svec{A}$ in every iteration. The intuitive expectation is, that methods from the first category have numerically cheap iterations, but superlinear convergence at best. The second category on the other hand, that includes all Newton--type solution strategies, requires the solution of several large linear systems but has the potential for quadratic convergence.

The PGJ, PGS and PSOR solvers can be seen as methods from the first category. Other methods have been suggested as well, such as projected gradient methods and spectral methods \cite{renouf2005,2013heynb}. The numerical results obtained with Nesterov's accelerated projected gradient descent method are very promising \cite{mazhar2015}.

Methods from the second category typically use so--called complementarity functions \cite{mangasarian1976,ferris1998,wanzke13}. A complementarity function is a function that is zero, if and only if the associated NCP is solved. This way, the complementarity problem is recast into the form of a root finding problem. Complementarity functions are typically nonsmooth and the gradient is not guaranteed to exist in general. Therefore the root finding problem must be solved using a semismooth Newton method, where in place of the gradients an arbitrary element from the subgradient of the complementarity function is taken.

In \cite{Daviet2011a} a hybrid strategy is proposed. The authors use a per--contact version of the Fischer--Burmeister complementarity function to model exact Coulomb friction. A semismooth Newton method per contact is used in an inner loop, while the global NCP is solved using a Gau\ss--Seidel--like outer loop.

An alternative strategy to overcome the nonsmoothness of the complementarity function is to approximate it using a smooth function. This is done by introducing a smoothing parameter $\alpha$ in such a way, that as $\alpha \to 0$, the smoothed function converges to the original complementarity function. Then a series of smooth root finding problems can be solved with a classical Newton method. Each subproblem is solved approximately, possibly with just one Newton iteration, before the smoothing parameter is decreased and the previous iterate serves as an initial guess for the next iteration. This strategy, though motivated completely differently, is very similar in spirit to interior point methods and path--following algorithms \cite{sun1999,chen2000}.

The use of Interior Point Methods (IPM) is proposed in \cite{Krabbenhoft2012,kleinert2014,kleinert2015,mangoni2018}. The strategy consists of using a logarithmic version of the objective function in the QP~\eqref{eq:qp} and adding a logarithmic penalty term, also called potential. The penalty term is infinitely large at the boundary (except at the exact solution) and drives the current iterate towards the interior of the feasible set. The zero--set of the potential is a smooth curve, called the central path. It can be parametrized using a parameter $\alpha$ in such a way, that the smooth curve passes through the solution of the complementarity problem at $\alpha = 0$. The strategy of interior point methods is to use a Newton--iteration to step towards the central path at a given $\alpha$ and then decrease the parameter iteratively. This way, the solution is approaches from within a close neighborhood of the central path. A big shortcoming of IPMs is, that the condition of the linear systems involved in the Newton--steps is known to diverge as $\alpha \to 0$. A remedy for this is preconditioning of the linear systems. Typically, a block--diagonal preconditioner is already very effective. IPMs are numerically costly, as a series of linear systems must be solved. But, especially in early iterations, the linear systems must not be solved to a high accuracy, so that inexact search directions can be calculated using Krylov subspace methods with generous tolerances. In \cite{corona2017} a compression based direct linear solver is proposed for the linear systems arising in IPMs in the context of nonsmooth rigid body dynamics.

Comparisons of some of the recent numerical methods for complementarity problems in the context of multibody simulations are provided in \cite{heyn2013,melanz2016,melanz2017,acary2018}.

\paragraph{Towards higher order integration}

The Moreau--Jean method discussed so far is attractive because it integrates nonsmooth motion at a fixed time step size for an arbitrary number of contacts, while maintaining robustness. Its main shortcoming is the global accuracy of order one, even in smooth phases of motion such as free flight without impacts. If no unilateral contacts are present in a multibody simulation, standard DAE solvers can be used with favorable properties such as second order accuracy or more, unconditional stability and controlled numerical damping.

Recently, some effort has been put into finding a compromise between both worlds. This is a mixed strategy that consists of a consistent time--stepping scheme to deal with unilateral contacts and a higher order integration method for the smooth phases of a simulation \cite{studer2008a,acary2012,chen2013,Schindler2014,schindler2015}. 

The general idea is to separate the time integration into two parts for smooth and nonsmooth motion. The smooth motion is integrated using a higher order DAE--solver on velocity level, i.e. using an index-2 formulation. If a unilateral constraint $g_i(\qe) \geq 0$ switches from inactive, $g_i(t_{k-1})>0$, to active, $g_i(t_k) \leq 0$, during the time step $[t_{k+1},t_k]$, the impact equations for this unilateral constraint are solved and the smooth motion is updated by the impulse contribution of the impact. 

During all time steps that contain a discontinuous velocity jump, the order cannot be expected to be higher than one. Therefore the global order of accuracy of the mixed time--stepping approach cannot exceed one. Locally, during time intervals without impacts however, the accuracy order is that of the DAE solver used for the smooth motion.

The above mentioned literature contains numerical experiments with the Newmark and Hugh--Hilbert--Taylor integrators, the generalized $\alpha$--method as well as several variants of the RADAU and Lobatto Runge--Kutta methods. In the context of higher order integration methods, constraint stabilization using the Gear--Gupta--Leimkuhler algorithm \cite{GeGL85}, see (\ref{ggl}), is discussed in \cite{schoeder2014,bruls2014,bruls2018}.

\subsubsection*{Acknowledgments}
We would like to thank Alessandro Tasora, Dan Negrut, 
and Klaus Dre{\ss}ler for fruitful discussions and valuable hints on 
constrained mechanical systems and nonsmooth extensions.

\bibliography{litbuch}
\bibliographystyle{jwbook}

\printindex

\end{document}